\documentclass[reqno,12pt]{amsart}

\address{
Department of Mathematics\\
University of British Columbia
}
\address{
Department of Mathematics\\
California Institute of Technology
}

\email{jbryan@math.ubc.ca, agholamp@its.caltech.edu}

\usepackage{amsmath,amsthm,amsfonts,amscd}
\usepackage{mathpazo}
\usepackage{verbatim}
\usepackage{afterpage}
\usepackage{xy}
\usepackage{graphicx}
\usepackage{enumerate}
\usepackage{diagrams}
\usepackage{amsmath,amsfonts,amssymb}
\usepackage{bbm,dsfont}
\usepackage{mathrsfs}
\usepackage{ifthen}
\usepackage{rotating}

\newarrow{Biject} <--->

\newtheorem{thm}{Theorem}

\newtheorem{lem}[thm]{Lemma}
\newtheorem{prop}[thm]{Proposition}
\newtheorem{conj}[thm]{Conjecture}

\newtheorem{lemma}[thm]{Lemma}

\theoremstyle{definition}
\newtheorem{rem}[thm]{Remark}

\newtheorem{defn}[thm]{Definition}
\newtheorem{remark}[thm]{Remark}

\newcommand{\Def}{\operatorname{Def}}
\newcommand{\age}{\operatorname{age}}
\newcommand{\LangRang}[1]{\left\langle #1 \right\rangle}
\renewcommand{\hat}{\widehat}
\newcommand{\M}{\overline{M}}
\newcommand{\X}{\mathcal{X}}
\newcommand{\mo}{\mathcal{O}_{\P^1}}
\newcommand{\Y}{\overline{Y}}
\newcommand{\Yh}{\widehat{Y}}
\newcommand{\Ys}{\widehat{Y}_{sing}}
\newcommand{\Ybs}{\widehat{\overline{Y}}_{\!\!sing}}
\newcommand{\Ybh}{\widehat{\overline{Y}}}
\newcommand{\Ybhg}{\widehat{\overline{Y}}_{\!\!gen}}
\newcommand{\Sh}{\widehat{S}}
\newcommand{\St}{\Sigma _{3}}
\newcommand{\Sth}{\widehat{\Sigma }_3}
\newcommand{\Ch}{\widehat{C}}
\newcommand{\pih}{\widehat{\pi}}
\newcommand{\qh}{\widehat{q}}

\newcommand{\Hilb}{\operatorname{Hilb}}
\newcommand{\Sym}{\operatorname{Sym}}
\newcommand{\dHilb}{\operatorname{-Hilb}}
\newcommand{\G}{G\dHilb (\C ^{3})}

\newcommand{\Gh}{\widehat{G}}

\newcommand{\Ghh}{\Gh\dHilb (\C ^{2})}
\newcommand{\g}{(g)}
\newcommand{\irr}{\operatorname{Irr}^*}
\newcommand{\Irr}{\operatorname{Irr}}
\newcommand{\con}{\operatorname{Conj}}
\newcommand{\CRC}{Crepant Resolution Conjecture}
\newcommand{\MC}{McKay correspondence}
\newcommand{\C}{\mathbb{C}}
\newcommand{\Q}{\mathbb{Q}}
\newcommand{\Z}{\mathbb{Z}}
\renewcommand{\P}{\mathbb{P}}
\newcommand{\GW}{Gromov-Witten}

\renewcommand{\ker}{\operatorname{Ker}}
\newcommand{\coker}{\operatorname{Coker}}

\newarrow{Into} {boldhook}--->
\newarrow{Corr} <--->

\newcommand{\epshmap}{
\begin{picture}(0,0)(-15,65)
\qbezier(10,-65)(0, -55)(0, 0) \qbezier(10,65)(0, 55)(0, 0)
\put(-10,0){\makebox(0,0)[t]{$\widehat{\epsilon }$}}
\put(11.5,-66){\vector(1,-1){0}}
\end{picture}
}

\newcommand{\epsbhmap}{
\begin{picture}(0,0)(-15,65)
\qbezier(10,-65)(0, -55)(0, 0) \qbezier(10,65)(0, 55)(0, 0)
\put(-10,0){\makebox(0,0)[t]{$\widehat{\overline{\epsilon}}$}} \put(11.5,-66){\vector(1,-1){0}}
\end{picture}
}

\newcommand{\pihmap}{
\begin{picture}(0,0)(13,63)
\qbezier(-17,60)(0, 40)(0, 0) \qbezier(-17,-63)(0,-20)(0, 0) \put(10,0){\makebox(0,0)[t]{$\pih $}}
\put(-19.6,-66){\vector(-1,-1){0}}
\end{picture}
}

\newcommand{\putlabelU}{
\begin{picture}(0,0)(0,0)
\put(5,1){\makebox(0,0)[t]{$ \scriptstyle U_1$}}
\end{picture}
}

\newcommand{\putlabelUU}{
\begin{picture}(0,0)(0,0)
\put(5,1){\makebox(0,0)[t]{$ \scriptstyle U_2$}}
\end{picture}
}

\newcommand{\putlabelUUU}{
\begin{picture}(0,0)(0,0)
\put(5,1){\makebox(0,0)[t]{$ \scriptstyle U_3$}}
\end{picture}
}

\newcommand{\putlabelV}{
\begin{picture}(0,0)(0,0)
\put(5,1){\makebox(0,0)[t]{$ \scriptstyle V_1$}}
\end{picture}
}\newcommand{\putlabelVV}{
\begin{picture}(0,0)(0,0)
\put(5,1){\makebox(0,0)[t]{$ \scriptstyle V_2$}}
\end{picture}
}\newcommand{\putlabelVVV}{
\begin{picture}(0,0)(0,0)
\put(5,1){\makebox(0,0)[t]{$ \scriptstyle V_1$}}
\end{picture}
}
\newcommand{\putlabelVVVV}{
\begin{picture}(0,0)(0,0)
\put(5,1){\makebox(0,0)[t]{$ \scriptstyle V_2$}}
\end{picture}
}
\begin{document}
\begin{abstract}
Let $G$ be a polyhedral group, namely a finite subgroup of
$SO(3)$. Nakamura's $G$-Hilbert scheme provides a preferred Calabi-Yau
resolution $Y$ of the polyhedral singularity $\C ^{3}/G$. The
classical {\MC} describes the \emph{classical} geometry of $Y$ in
terms of the representation theory of $G$. In this paper we describe
the \emph{quantum} geometry of $Y$ in terms of $R$, an ADE root system
associated to $G$.  Namely, we give an explicit formula for the
Gromov-Witten partition function of $Y$ as a product over the positive
roots of $R$. In terms of counts of BPS states (Gopakumar-Vafa
invariants), our result can be stated as a correspondence: each
positive root of $R$ corresponds to one half of a genus zero BPS
state. As an application, we use the {\CRC} to provide a full
prediction for the orbifold {\GW} invariants of $[\C^3/G]$.
\end{abstract}

\title{The Quantum McKay Correspondence for polyhedral singularities}
\author{Jim Bryan and Amin Gholampour} \maketitle
\markright{The Quantum McKay Correspondence}


\section{Introduction}\label{sec: intro}

Let $G$ be a finite subgroup of $SU (3)$, and let $X $ be the quotient
singularity
\[
X=\C ^{3}/G.
\]
There is a preferred Calabi-Yau resolution 
\[
\pi :Y\to X
\]
given by Nakamura's
$G$-Hilbert scheme (Definition~\ref{defn: G-Hilb})
\[
Y=\G .
\]
The classical McKay correspondence describes the geometry of $Y$ in
terms of the representation theory of $G$ \cite{BKR, McKay}. One of
the original formulations\footnote{The modern way to formulate the
correspondence is an equivalence of derived categories \cite{BKR}.} of
the correspondence is a bijection of finite sets:

\bigskip
\noindent 
\textbf{Classical McKay Correspondence:}
\begin{equation*}
\{\text{Irreducible representations of $G$} \} \rBiject \text{Geometric basis for $H^{*} (Y,\Z )$}.
\end{equation*}
\smallskip

We consider the case where $G$ is a \emph{polyhedral group}, namely a
finite subgroup of $SO (3)\subset SU (3)$. Such groups are classified
into three families $A_{n}$, $D_{n}$, and $E_{n}$, given by the cyclic
groups, the dihedral groups, and the symmetries of the platonic
solids. In this paper we describe the \emph{quantum} geometry of $Y$,
namely its Gromov-Witten theory, in terms of the associated ADE root
system $R$. Our main theorem provides a closed formula for the
Gromov-Witten partition function in terms of the root system $R$
(Theorem~\ref{thm:main theorem}).

This result can also be formulated in terms of the Gopakumar-Vafa
invariants which are defined in terms of Gromov-Witten invariants and
physically correspond to counts of BPS states. Our result provides a
\emph{quantum McKay correspondence} which can be stated as a
natural bijection of finite sets:

\bigskip
\noindent 
\textbf{Quantum McKay Correspondence:}
\[
\left\{\begin{array}{c}
\text{Positive roots of $R$}\\
\text{other than binary roots}
\end{array} \right\}
\lBiject
\left\{\begin{array}{c}
\text{Contributions of $\frac{1}{2}$ to}\\
\text{BPS state counts of $Y$}
\end{array} \right\}
\]

\bigskip

The \emph{binary roots} are a distinguished subset of the simple roots
and are given in Figure~\ref{fig:ADE}.  In the above correspondence,
all the BPS state counts are in genus 0; we show that $Y$ has no
higher genus BPS states. The meaning of the above bijection is that
$n^{0}_{\beta } (Y)$, the Gopakumar-Vafa invariant which counts genus
zero BPS states in the curve class $\beta $, is given by half the
number of positive roots corresponding to the class $\beta $. We will
explain below how the positive roots of $R$ (other than binary roots)
naturally correspond to curve classes in $Y$. In general, many
positive roots may correspond to the same curve class and moreover,
that curve class may be reducible and/or non-reduced. Thus the above
simple correspondence leads to a complicated spectrum of BPS states:
they can be carried by complicated curves and have various
multiplicities. As an example, we work out the $D_{5}$ case explicitly
in \S\ref{sec: D5 case}.

Our method for computing the Gromov-Witten partition function of $Y$
uses a combination of localization, degeneration, and deformation
techniques. There is a singular fibration $Y\to \C $ by (non-compact)
$K3$ surfaces and we use localization and degeneration to relate the
partition function of $Y$ to the partition function of a certain
smooth $K3$ fibration $\Yh\to \C $.  We compute the partition function
of $\Yh $ by deformation methods. The root system $R$ makes its
appearance in the versal deformation space of the central fiber of
$\Yh \to \C$. We give a brief overview of the proof in \S~\ref{subsec: overview of proof}.

\subsection{The main result.}

Since $G$ is a finite subgroup of $SO (3)$, it preserves the quadratic
form $t=x^{2}+y^{2}+z^{2}$ and hence $G$ acts fiberwise on the
family of surfaces 
\begin{align*}
\C^3 &\rightarrow \C \quad\\
(x,y,z) &\mapsto x^2+y^2+z^2.
\end{align*}

Let
\[
Q_t=\{x^2+y^2+z^2=t\} \subset \C^3
\]
be the fiber over $t \in \C$.  The central fiber is a quadric cone,
which is isomorphic to the $A_{1}$ surface singularity:
\[
Q_0 \cong \C^2/\{\pm 1\}.
\]

Via the functorial properties of the $G$-Hilbert scheme, $Y$ inherits a map
\[
\epsilon: Y \rightarrow \C
\]
whose general fiber, $S_t$, is the minimal resolution of $Q_t/G$, and
whose central fiber, $S_0$, is a partial resolution of
$$Q_0/G \cong \C^2/\Gh,$$ where $\Gh$ is the binary version of $G$, namely the
preimage of $G$ under the double covering $j$:
\[
\begin{diagram}[height=1.5em,width=1.5em]
\Gh&\rInto & SU (2)\\
\dTo & &\dTo_j \\
G&\rInto & SO (3)
\end{diagram}
\]

Let $\Sh =\Ghh $ be the minimal resolution of the surface singularity
$\C ^{2}/\hat{G}$
\[
\pih: \Sh  \rightarrow \C^2/\Gh.
\]
The map $\hat{\pi }$ factors through $\Sh \to S_{0}$ which gives rise
to a map
\[
f:\Sh \to Y.
\]
See Figure~\ref{fig:big diagram} in \S\ref{sec:proof of the main thm}
for a diagram of these maps.

The map $f$ is constructed modularly by Boissiere and Sarti
\cite{Boissiere-Sarti} who prove a beautiful compatibility between the
McKay correspondences for $\hat{\pi }:\Sh \to \C ^{2}/\hat{G}$ and
$\pi :Y\to \C ^{3}/G$. The compatibility can be expressed as the
commutativity of the following diagram of set maps.  \medskip
\[
\begin{diagram}
\left\{\begin{array}{l}
\text{Non-trivial irreducible}\\
\text{$\hat{G}$-representations, $\hat{\rho }$}
\end{array} \right\}
&\rBiject&
\left\{\begin{array}{l}
\text{Components}\\
\text{of $\hat{\pi }^{-1} (0)$,  $\hat{C}_{\hat{\rho }}$ }
\end{array} \right\}\\
\uTo_{\text{pullback}}&&\uTo_{
\begin{array}{l}
\text{proper transform}\\
\text{under $f$}
\end{array}
}\\
\left\{\begin{array}{l}
\text{Non-trivial irreducible}\\
\text{$G$-representations, $\rho $}
\end{array} \right\}
&\rBiject&
\left\{\begin{array}{l}
\text{Components}\\
\text{of ${\pi }^{-1} (0)$, $C_{\rho }$}
\end{array} \right\}\\
\end{diagram}
\]

\medskip The vertical arrows are injections and horizontal arrows
are the bijections from the McKay correspondences in dimensions two
and three. The bijections say that the irreducible components of the
exceptional fibers over 0 (which are smooth rational curves in these
cases) are in natural correspondence with non-trivial irreducible
representations of $\hat{G}$ and $G$ respectively.  The commutativity
of the diagram says that the map $f$ contracts exactly those
exceptional curves which correspond to representations of $SU (2)$ not
coming from representations of $SO (3)$.



The configurations of curves in $\hat{\pi
}^{-1} (0)$ and $\pi ^{-1} (0)$ are given in
\cite[Figures~5.1\&~5.2]{Boissiere-Sarti}.  The curves $\hat{C}_{\hat{\rho } }$
are smooth rational curves meeting transversely with an intersection
graph of ADE type. This gives rise to an ADE root system $R$ and a
natural bijection between the simple roots of $R$ and the irreducible
components of $\hat{\pi }^{-1} (0)$.
\[
\begin{diagram}
\left\{\begin{array}{l}
\text{Non-trivial irreducible}\\
\text{$\hat{G}$-representations, $\hat{\rho }$}
\end{array} \right\}
&\rBiject&
\left\{\begin{array}{l}
\text{Components}\\
\text{of $\hat{\pi }^{-1} (0)$, $\hat{C}_{\hat{\rho }}$}
\end{array} \right\}
&\rBiject &
\left\{\begin{array}{l}
\text{Simple roots}\\
\text{of $R$, $e_{\hat{\rho }}$}
\end{array} \right\}
\end{diagram}
\]

\begin{figure} 
\centering \includegraphics{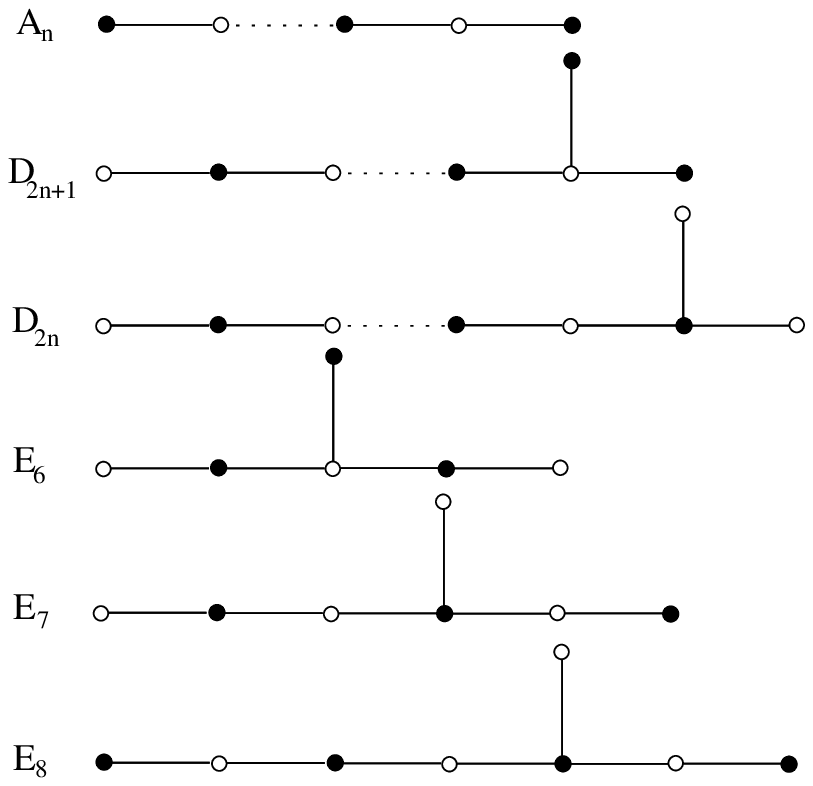} \caption{The set of \emph{binary
roots} of an ADE root system is the subset of the simple roots
corresponding to the black nodes in the above Dynkin diagrams. The ADE
Dynkin diagrams are dual to the configurations of exceptional curves
in $\Sh $. The map $f:\Sh \to Y$ contracts exactly the curves
corresponding to binary roots.}\label{fig:ADE}
\end{figure}

For any positive root 
\[
\alpha=\sum _{\hat{\rho }}\alpha ^{\hat{\rho }}e_{\hat{\rho }} \in R^{+}
\]
we define maps
\[
\hat{c}:R^{+}\to H_{2} (\hat{S},\Z )\quad \quad \quad c:R^{+}\to H_{2} (Y,\Z )
\]
by
\[
\hat{c} \left(\alpha \right) = \sum _{\hat{\rho }}\alpha ^{\hat{\rho
}} \hat{C}_{\hat{\rho }}\quad \quad \quad 
c (\alpha ) = f_{*} (\hat{c} (\alpha )).
\]
Note that $\hat{c}$ is injective but $c$ is not since $f_{*}$
contracts curves. By definition, a simple root is \emph{binary} if it
is in the kernal of $c$, or in other words, if the corresponding curve
is contracted by $f:\Sh \to Y$. See Figure~\ref{fig:ADE} for a case by
case description of binary roots.

Let $N_{\beta }^{g} (Y)$ be the genus $g$ Gromov-Witten invariant of
$Y$ in the class $\beta \in H_{2} (Y,\Z )$. It is defined to be the
degree of the virtual fundamental cycle on $\M_{g} (Y,\beta )$, the
moduli space of genus $g$ stable maps to $Y$ in the class of $\beta $
\cite{coxkatz,mirrorbook}. In the event that $\M _{g} (Y,\beta )$ is
non-compact, as it is for some classes $\beta $, the invariant
$N_{\beta }^{g} (Y)$ is defined via virtual localization
\cite[\S2]{Br-Pa-local-curves} using the natural $\C ^{*}$ action on
$Y$. The action on $Y$ is induced by the diagonal action of $\C ^{*}$
on $\C ^{3}$.

The Gromov-Witten invariants of non-zero degree are assembled into a
generating function called the \emph{reduced Gromov-Witten partition
function}
\[
Z^{Y} (q,\lambda ) = \exp\left(\sum _{g=0}^{\infty }\sum _{\beta \neq
0} N_{\beta }^{g} (Y) q^{\beta } \lambda ^{2g-2} \right).
\]

The main theorem of this paper can then be written as follows.
\begin{thm} \label{thm:main theorem} Let $G$ be a finite subgroup of
$SO(3)$. The reduced {\GW} partition function of $Y=\G$ is given by
\[
Z^Y(q,\lambda)=\prod_{\begin{smallmatrix}\alpha \in R^{+}\\
c (\alpha )\neq 0  \end{smallmatrix}} \prod_{m=1}^{\infty}\left(1-q^{c (\alpha )} \left( -e^{i\lambda }\right)^{m} \right)^{m/2}.
\]
\end{thm}

\subsection{Applications.} There are very few Calabi-Yau threefolds
for which the Gromov-Witten partition function has been completely
determined. In addition to the $G$-Hilbert schemes of this paper, the
only other cases are toric Calabi-Yau threefolds and local curves. The
partition function of a toric Calabi-Yau threefold can be computed via
the topological vertex formalism \cite{AKMV,Li-Liu-Liu-Zhou} and the
partition function of a local curve is determined in
\cite{Br-Pa-local-curves}. Having an explicit exact formula for the
Gromov-Witten partition function leads to a number of applications
some of which we describe below.

\subsubsection{Don\-ald\-son-Thomas and Pand\-hari\-pan\-de-Thomas
theory} The Don\-ald\-son-Thomas invariants of a threefold are defined
in terms of the moduli space of ideal sheaves on the threefold. It was
conjectured in \cite[Conj.~3]{MNOP1} and
\cite[Conj.~3R]{Br-Pa-local-curves} that the reduced partition
function of {\GW} and Don\-ald\-son-Thomas theory are equal after the
change of variables $-e^{i\lambda}=Q$. Recently, a conjecturally
equivalent theory using a moduli space of stable pairs has been
defined by Pand\-hari\-pan\-de and Thomas
\cite{Pandharipande-Thomas1}. Since the dependence of the genus
parameter $\lambda $ in our formula for the Gromov-Witten partition
function of $Y$ is explicitly through the quantity $-e^{i\lambda }$,
Theorem \ref{thm:main theorem} immediately yields a prediction for the
reduced Don\-ald\-son-Thomas partition function of $Y$ and the
Pand\-hari\-pan\-de-Thomas partition function of $Y$.

\subsubsection{BPS states} The result of Theorem \ref{thm:main
theorem} in terms of the BPS state counts is remarkably simple. BPS
state counts, defined physically by Gopakumar and Vafa \cite{Go-Va},
can be defined mathematically in terms of Gromov-Witten invariants
\cite[Def.~1.1]{Br-Pa}. Namely, the genus $g$, degree $\beta$
Gopakumar-Vafa invariants of $Y$, denoted by $n^g_\beta (Y)$, are
defined by the formula
\[
\sum_{g=0}^\infty \sum_{\beta \neq 0} N^g_\beta(Y) q^\beta
\lambda^{2g-2}=\sum_{g=0}^\infty \sum_{\beta \neq 0} n^g_\beta (Y)
\sum_{d>0}
\frac{1}{d}\left(2\sin\left(\frac{d\lambda}{2}\right)\right)^{2g-2}q^{d\beta}.
\]

It is well known that each genus zero BPS state $n^{0}_{\beta }$
contributes a factor of 
\[
\prod _{m=1}^{\infty }\left(1-q^{\beta } (-e^{i\lambda })^{m}
\right)^{m n_{\beta }^{0}}
\]
to the partition function (e.g. \cite[Proof of
Thm~3.1]{Behrend-Bryan}) and so Theorem \ref{thm:main theorem} implies
that BPS states of $\G$ are given by
\[
n^g_\beta (Y) =\begin{cases}  \frac{1}{2}|c^{-1}(\beta)| & \quad g=0, \\
0 & \quad g>0, \end{cases}
\]
where $|\cdot |$ denotes the cardinality of a finite set.  This makes
precise the Quantum McKay Correspondence which we described
earlier. Namely, each non-binary positive root $\alpha $ contributes
$\frac{1}{2}$ to the genus 0 BPS state count in the class $\beta =c
(\alpha )$. Note that the binary roots are precisely those roots
$\alpha $ with $c (\alpha )=0$. The example of $D_{5}$ is worked out
explicitly in \S~\ref{sec: D5 case}.

\begin{rem}
It can be directly observed that for any curve class $\beta$ in $Y$,
$|c^{-1}(\beta)|\in \{0,1,2,4,8 \}$. In particular, for some of the
curve classes the BPS state count is $1/2$. This does not contradict
the conjectural integrality of the BPS state counts of the Calabi-Yau
threefolds. This is because of the non compactness of $Y$ that causes
the ordinary {\GW} invariant corresponding to some particular curve
classes not to be a priori well defined. We use the $\C^*$ action on
$Y$ to define the equivariant version of the invariants for such curve
classes. The equivariant {\GW} invariants agree with the ordinary
{\GW} invariants whenever the latter are well-defined.
\end{rem}

\subsubsection{Orbifold {\GW} theory of $[\C^3/G]$}
\label{sec:orbifold GW theory} Let $\X$ be the orbifold associated to
the singular space $\C^3/G$.  In \S\ref{sec:CRC} we show that Theorem
\ref{thm:main theorem}, together with classical invariants of $Y$ and
the {\CRC} \cite{Bryan-Graber,Coates-Ruan}, give a prediction for the
genus zero {\GW} invariants of $\X$. To express this, we introduce
some notation.

The inertia orbifold, $I\X$ is a union of contractible connected
components, $I\X_{\g}$ indexed by $\g \in \con(G)$, the set of the
conjugacy classes of $G$. Recall that the (equivariant) orbifold
cohomology of $\X$ is by definition \cite{Chen-Ruan}
\begin{align*}
H^*_{orb}(\X)=\bigoplus_{(g) \in\; \con(G)}
H_{\C^*}^{*-2\imath_{\g}}(I\X_{\g}),\end{align*} where $\imath_{\g}$
is the degree shifting number (age) of the component indexed by $\g
\in \con(G)$. 

It is a consequence of $G\subset SO (3)$ that all non-trivial elements
of $G$ act on $\C ^{3}$ with age one (c.f. Lemma~\ref{lem: hard
lefschetz => G in SO(3) or SU(2)}). Therefore, there exists a
canonical basis $\{\delta_{(g)}\}_{(g) \in\; \con(G)},$ for
$H^*_{orb}(\X)$ where
\[
\delta_{(e)}\in H^0_{orb}(\X) \quad \text{ and
} \quad \delta_{(g)}\in H^2_{orb}(\X) \text{ for } \g \neq (e).
\]
Let $x=\{x_{(g)}\}_{(g) \in\; \con(G)}$ be a set of variables
parameterizing this cohomology basis. For any given vector
$n=(n_{\g})_{\g \in\; \con(G)}$ of nonnegative integers, we use the
following notation $$\delta^n= \prod_{\g \in\; \con(G)}
\delta_{\g}^{n_{\g}} \quad \text{and} \quad \frac{x^n}{n!}=\prod_{\g
\in\; \con(G)} \frac{x_{\g}^{n_{\g}}}{n_{\g}!}.$$ Suppose
\[
\displaystyle |n|=\sum_{\g \in\; \con(G)} n_{\g}.
\]
We denote the
genus zero, $|n|$-point, equivariant orbifold {\GW} invariants of
$\X$ corresponding to the vector $n$ by $\langle \delta^n \rangle$
(see \cite{Chen-Ruan,Ab-Gr-Vi} and \cite[\S1.4]{Bryan-Graber}). 

\begin{defn} Using the notation above, we write the genus zero orbifold
{\GW} potential function of $\X$ as
$$F^\X(x)=\sum_{n\in \Z_+^r} \langle \delta^n  \rangle \, \frac{x^n}{n
!}$$ where $r=|\con(G)|$. For simplicity, we set the variable
$x_{(e)}$ to zero (see Remark \ref{rem:x0 part}).
\end{defn}

In \S\ref{sec:CRC}, we establish the following prediction for
$F^\X(x)$. The reader can find a similar prediction for the orbifold
{\GW} invariants of $[\C^2/\Gh]$ given in
\cite[Conjecture~11]{Bryan-Gholampour2}.
\begin{conj} \label{conj:orbifold prediction}
Let $G$ be a finite subgroup of $SO(3)$ and let $\mathbf{h} (s)$ be a
series defined by
\[
\mathbf{h}'''(s)= \frac{1}{2}\tan\left(-\frac{s}{2}\right).
\]
Let
\[
\mathsf{X}_{\rho } = \frac{1}{|G|}\left(2\pi \dim (\rho )+
\sum_{g \in\, G}\sqrt{3-\chi_V (g)}\;\chi_\rho (g) x_{\g} \right)
\] 
where $\chi _{\rho } (g)$ is the character of the $G$-representation
$\rho $ evaluated on $g$ and $V$ is the three dimensional
representation induced by the embedding $G\subset SO (3)\subset SU
(3)$.

Then the genus zero orbifold {\GW}
potential function of $\X=[\C^3/G]$ is given by
\[
F^{\X}(x)=\frac{1}{2}\sum_{\alpha \in\, R^+} \mathbf{h} \left(\pi +
\sum _{\rho}
\alpha ^{\rho }\,\mathsf{X}_{\rho }
\right)
\]
where the second sum is over non-trivial irreducible
$G$-representations $\rho $ and $\alpha ^{\rho }$ is the coefficient
of $e_{\rho }$ (the simple root corresponding to $\rho $) in the
positive root $\alpha $.
\end{conj}
Note that terms of order less than three in $\mathbf{h} (s)$, and
consequently in $F^{\X } (x)$, are not defined. This reflects the fact
that the invariants with fewer than three insertions are also not
defined. The orbifold Gromov-Witten invariants of $[\C ^{3}/G]$ may be
expressed in terms of $G$-Hurwitz-Hodge integrals. Hurwitz-Hodge
integrals have been the topic of extensive recent research
\cite{Bryan-Gholampour1, Bryan-Graber-Pandharipande,
Cavalieri-Hurwitz-Hodge, Chiodo, CCIT-CRC,
Johnson-Panharipande-Tseng,Zhou-HH}.

In \S\ref{sec:CRC} we prove the following proposition.

\begin{prop}\label{prop: Fx conjecture implies CRC}
If Conjecture~\ref{conj:orbifold prediction} is true, then the genus
zero crepant resolution conjecture holds for $Y\to \C ^{3}/G$.
\end{prop}

Conjecture~\ref{conj:orbifold prediction} has been proved for $G$
equal to $\Z _{2}\times \Z _{2}$ or $A_{4}$ in
\cite{Bryan-Gholampour1} (which correspond to $D_{4}$ and $E_{6}$ in
the ADE classification), and for $\Z _{n+1}$ in \cite{CCIT-CRC} (which
corresponds to $A_{n}$ in the ADE classification). Zhou has recently
proved a higher genus analog of Conjecture~\ref{conj:orbifold
prediction} in the $A_{n}$ case \cite{Zhou-An}.

\section{Proof of Theorem~\ref{thm:main theorem}.}\label{sec:proof of
the main thm}

\subsection{Overview}\label{subsec: overview of proof} To compute the
Gromov-Witten partition function of $Y$, we construct and study four
families of (non-compact) $K3$ surfaces:
\begin{align*}
\epsilon :Y&\to \C \quad \quad \, \hat{\epsilon }:\hat{Y}\to \C \\
\overline{\epsilon} :\overline{Y}&\to \P ^{1}\quad \quad \hat{\overline{\epsilon }}:\hat{\overline{Y}}\to \P ^{1} 
\end{align*}
We obtain the smooth family $\hat{\epsilon }$ from $\epsilon $ by
taking the double cover of $Y$ branched along the (singular) central
fiber and then taking a small resolution of the resulting conifold
singularities. The family $\overline{\epsilon }$ extends the family
$\epsilon $ to the compact base $\P^1$. The family
$\hat{\overline{\epsilon }}$ extends the family $\hat{\epsilon }$ to
the compact base $\P^1$.

We study $Z^{Y} (q,\lambda )$ and $Z^{\hat{Y}} (\hat{q},\lambda )$,
the reduced Gromov-Witten partition functions of $Y$ and $\hat{Y}$
respectively. We also study $Z^{\overline{Y}}_{f} (q,\lambda )$ and
$Z^{\hat{\overline{Y}}}_{f} (\hat{q},\lambda )$, the fiber class
reduced Gromov-Witten partition functions of $\overline{Y}$ and
$\hat{\overline{Y}}$ respectively.

We denote the set of non-trivial irreducible representations of $G$
and $\hat{G}$ by $\irr(G)$ and $\irr(\Gh)$, respectively. Using the
identifications in Lemma~\ref{lem: homology isos}, the quantum
variables of $Y$ and the fiber class quantum variables of
$\overline{Y}$ are given by $\{q_{\rho } \}_{\rho \in \Irr^{*} (G)}$
whereas the quantum variables of $\hat{Y}$ and the fiber class quantum
variables of $\hat{\overline{Y}}$ are given by $\{\hat{q}_{\hat{\rho
}} \}_{\hat{\rho }\in \Irr ^{*}(\hat{G})}$.

There is a specialization of the variables given by
\[
\hat{q}_{\hat{\rho}} = \begin{cases}
q_{\rho } &\text{if $\,\hat{\rho }\,$ pulls back from $\rho \in \Irr ^{*} (G)$}\\
1&\text{if $\,\hat{\rho }\,$ is not the pullback of a $G$ representation.}
\end{cases}
\]
We will prove the following four equalities
\begin{align}
Z^{{Y}} ({q},\lambda )^{2} & = Z_{f}^{{\overline{Y}}} ({q},\lambda )\label{eqn: ZY2=ZYbar}\\
Z_{f}^{{\overline{Y}}} ({q},\lambda
)^{2}&=\frac{Z_{f}^{\hat{\overline{Y}}}
(\hat{q},\lambda )}{ Z_{f}^{\Ybh } (\hat{q},\lambda )|_{q=0}}  \label{eqn: ZYbar2=ZYbarhat}\\
 Z_{f}^{\hat{\overline{Y}}} (\hat{q},\lambda )\,\,& =Z^{\hat{Y}} (\hat{q},\lambda )^{2} \label{eqn: ZYhat2=ZYbarhat}\\
Z^{\hat{Y}} (\hat{q},\lambda )\,\, & = 
\prod _{
\alpha \in R^{+}
}\prod _{m=1}^{\infty }
\left(1-\hat{q}^{\,\,\hat{c} (\alpha ) }\left(-e^{i\lambda } \right)^{m}
\right)^{m}\label{eqn: ZYhat}
\end{align}
where the $\hat{q}\mapsto q$ specialization is implicit in
equation~\eqref{eqn: ZYbar2=ZYbarhat}. The main theorem follows
immediately from these four equations.

In \S\ref{subsec: Geometry of Y and the fibrations} we discuss the
geometry of the $G $-Hilbert scheme $Y$ and we construct the four
families. In Proposition~\ref{prop: formula for Z of Yhat} we prove
\eqref{eqn: ZYhat} using a deformation argument. In Lemma~\ref{lem:
relating Z of families over C and P1} we prove \eqref{eqn:
ZYhat2=ZYbarhat} and \eqref{eqn: ZY2=ZYbar} using a localization
argument. In Proposition~\ref{prop:relation to conifold resolution} we
prove \eqref{eqn: ZYbar2=ZYbarhat} using a degeneration argument and
the Li-Ruan formula for the Gromov-Witten invariants of a conifold
transition \cite{Li-Ruan}.

\subsection{Geometry of $Y=\G$ and the families $\epsilon $,
$\hat{\epsilon }$, $\overline{\epsilon }$, $\hat{\overline{\epsilon
}}$} \label{subsec: Geometry of Y and the fibrations} ~

We start by recalling the definition of $Y$, our main object of study.
\begin{defn}\label{defn: G-Hilb}
Let $Y=\G$ be Nakamura's $G$-Hilbert scheme \cite{Nakamura}, namely
the subscheme of $\Hilb^{|G|}(\C^3)$ parameterizing closed
$G$-in\-var\-i\-ant subschemes $Z \subset \C^3$ such that
$H^0(\mathcal{O}_{Z})$ is isomorphic to the regular representation of
$G$. Here $\Hilb^{|G|}(\C^3)$ denotes the Hilbert scheme of length
$|G|$ subschemes of $\C^3$.
\end{defn} 

Consider the Hilbert-Chow morphism from $\Hilb^{|G|}(\C^3)$ to
the $|G|$-th symmetric product of $\C^3$
\[
h: \Hilb^{|G|}(\C^3) \rightarrow \Sym^{|G|}(\C^3)
\]
that sends a closed subscheme $Z \subset \C^3$ to the zero cycle $[Z]$
associated to it. The morphism $h$ is projective and it restricts to a
projective morphism
\[
\pi: Y \rightarrow \C^3/G.
\]
Bridgeland, King, and Reid prove that $Y$ is smooth, reduced and
irreducible variety, and that $\pi$ is a crepant resolution
\cite[Theorem~1.2]{BKR}.

\begin{rem} \label{rem:H-Hilb}
More generally, one can define the $H$-Hilbert scheme
\mbox{$H$-Hilb$(M)$}, where $M$ is any smooth $n$-dimensional
quasi-projective variety, and $H$ is a finite group of automorphisms
of $M$, such that the canonical bundle of $M$ is a locally trivial
$H$-sheaf \cite{Nakamura,Reid-asterisque,BKR}.
\begin{enumerate} [(i)]
\item \mbox{$H$-Hilb$(M)$} represents the moduli functor of
$H$-clusters \cite[\S3]{Boissiere-Sarti}.
\item If $n \le 3$ then $H$-$\Hilb (M)\to M/G$ is a crepant
resolution \cite[Theorem~1.2]{BKR}.
\end{enumerate}
\end{rem}

\begin{center}
\begin{figure}[ht]
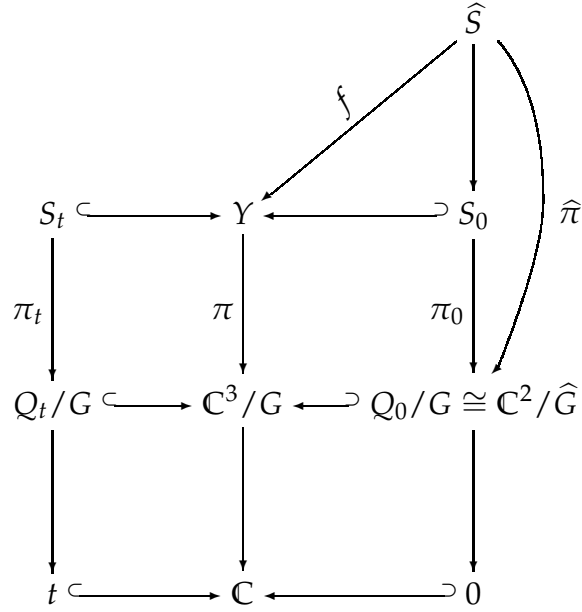

\begin{diagram}
\\   &        &    &   & \Sh& \pihmap\\
     &        &      & \ldTo^{f} & \dTo&\\
S_t  & \rInto & Y  &\lInto & S_0\\
\dTo^{\pi_t} &  & \dTo^\pi & & \dTo^{\pi_0}& \\
Q_t/G &  \rInto & \C^3/G &\lInto & Q_0/G\cong \C^2/\Gh& \\
\dTo &  & \dTo & & \dTo&\\
 t & \rInto & \C & \lInto & 0& \\
\end{diagram}
\caption{The family of surfaces $\epsilon :Y\to \C $ and related morphisms.} \label{fig:big diagram}
\end{figure}
\end{center}

Recall that $\epsilon :Y\to \C $ is the composition of the map $\pi
:Y\to \C ^{3}/G$ with the map $\C ^{3}/G\to \C$ given by
$(x,y,z)\mapsto x^{2}+y^{2}+z^{2}$. We analyze the fibers
$S_{t}=\epsilon ^{-1} (t)$. The relevant morphisms are shown in the
diagram in Figure~\ref{fig:big diagram}.

\afterpage{}

By
\cite[Theorem~1.1]{Boissiere-Sarti},
\[
\pi: Y \rightarrow \C^3/G
\]
restricts to a partial resolution,
\[
\pi_0:S_0 \rightarrow Q_0/G,
\]
whose exceptional locus coincides with $\pi^{-1}(0) $. 

Let $\Sh =\Ghh $.
Boissiere and Sarti \cite[Theorems~1.1\&
8.1]{Boissiere-Sarti} relate the geometry of $\Sh$ and $Y$ by constructing a morphism
\[
f:\Sh \rightarrow Y
\]
and proving it factors through $S_0 \subset Y$ as the minimal
resolution. They prove it maps $\Ch_\rho$ isomorphically to $C_\rho$
if $\rho \in \irr(G)$, and it contracts $\Ch_\rho$ if $\rho \in
\irr(\Gh) \backslash \irr(G)$. By the functoriality of the $G$-Hilbert
scheme (see Remark \ref{rem:H-Hilb} (i)), for each $t\neq 0$ $$S_t
\cong \mbox{$G$-Hilb$(Q_t)$},$$ and by Remark \ref{rem:H-Hilb} (ii),
$\pi$ restricts to the minimal resolution
\[
\pi_t:S_t \rightarrow Q_t/G.
\]

\begin{rem} \label{rem:C* on eps}
The action of $s \in \C ^{*}$ on $\C ^{3}$ is given by
\[
(x,y,z)\mapsto (s x,s y,s z).
\]
Thus in order to make the map $t=x^{2}+y^{2}+z^{2}$ equivariant, the
action of $s $ on the base of the family $\epsilon $ is given by
\[
t\mapsto s^{2}t.
\]
\end{rem}

The family $\epsilon :Y\to \C $ has a singular central fiber
$S_{0}$. We construct a related smooth family $\hat{\epsilon }$ as
follows.  Let
\[
\Ys \rightarrow \C
\]
be the double cover of $Y$ branched over the central fiber $S_0$. The
singularities of $S_{0}$ are obtained by contracting $-2$ curves in
$\hat{S}$, which are disjoint by
\cite[Figures~5.1\&5.2]{Boissiere-Sarti} (c.f.
Figure~\ref{fig:ADE}). Hence the singularities of $S_{0}$ are isolated
rational double point singularities of type $A_{1}$. Consequently,
$\Ys$ has conifold singularities at the singular points of $S_0$, and
it is smooth away from these points. Let $\Yh \rightarrow \Ys$ be a
small resolution of the conifold singularities. Thus we obtain a
$\C^*$-equivariant family
\[
\widehat{\epsilon}:\Yh \rightarrow \C
\]
that makes the following diagram of $\C ^{*}$-equivariant maps
commute.
\begin{diagram}
\epshmap &\Yh   &  & \\
&\dTo & \rdTo^{\hat{f}} & \\
&\Ys  & \rTo & Y \\
&\dTo &  & \dTo_\epsilon \\
&\C &  \rTo^{t\mapsto t^{2}} & \C\\
\end{diagram}

By the uniqueness of the minimal resolution in dimension two, the
central fiber $\widehat{\epsilon}^{\;-1}(0)$ is isomorphic to $\Sh$,
the minimal resolution of $\C^2/\Gh$. Hence, $\widehat{\epsilon}$ is a
$K3$ fibration with all the fibers being smooth, and $\{\Ch_{\hat{\rho
}}\}_{\hat{\rho } \in \irr(\Gh)}$ forms a basis for
\[
H_2(\Yh,\Z)\cong H_2(\Sh,\Z) 
\]
(see Lemma~\ref{lem: homology isos}).
\begin{rem} \label{rem:C* on epshat}
Note that because the base of the family $\widehat{\epsilon}$ double
covers the base of $\epsilon$, the $\C^*$-action has weight one on the
base of $\widehat{\epsilon}$.
\end{rem}

We construct $\overline{\epsilon }$ and $\hat{\overline{\epsilon }}$
by compactifying the base of $\epsilon $ and $\hat{\epsilon }$ as
follows.  Let
\[
W=\C^3 \bigcup_{\C^3\backslash Q_0}
\C^3
\]
where the two copies of $\C^3$ are glued along
$\C^3\backslash Q_0$ via the map
\begin{equation}\label{eqn: x-->x/(x^2+y^2+z^2)}
(x,y,z) \mapsto
\left(\frac{x}{x^2+y^2+z^2},\frac{y}{x^2+y^2+z^2},\frac{z}{x^2+y^2+z^2}\right).
\end{equation}
The two copies of the map $\varepsilon: \C^3 \rightarrow \C$ given by
$(x,y,z)\mapsto x^{2}+y^{2}+z^{2} $ patch together to give the
map 
\[
\overline{\varepsilon}:W \rightarrow \mathbb{P}^1.
\]
The actions of $\C^*$ and $G$ extends naturally to
$\overline{\varepsilon}$. We define 
\[
\overline{Y}=G\dHilb (W).
\]
This construction gives rise to the $\C ^{*}$--equivariant family
\[
\overline{\epsilon}:\Y \rightarrow \P^{1}.
\]
Restricted to each of the $\C^*$--invariant affine patches of $\P
^{1}$, the family is isomorphic to $\epsilon $, but with opposite $\C
^{*}$ weights.

We can easily extend the construction of $\hat{\epsilon }$ to the
compactified base as follows.

Let 
\[
\Ybs \rightarrow \mathbb{P}^1
\]
be the double cover of the
family $\overline{\epsilon}$ branched over the the fibers
$\overline{\epsilon}^{\;-1}(0)$ and
$\overline{\epsilon}^{\;-1}(\infty)$. Let 
\[
 \Ybh \rightarrow \Ybs
\]
be a small resolution of the conifold singularities (chosen as stated
in Remark~\ref{rem: Ybarhat constructed by gluing} below). One obtains
the family
\[
\widehat{\overline{\epsilon}}:\Ybh \rightarrow \mathbb{P}^1
\]
of smooth non-compact $K3$ surfaces that make the following
diagram of $\C ^{*}$-equivariant maps commute.
\begin{diagram}
\epsbhmap &\Ybh   &  & \\
&\dTo & \rdTo^{\hat{\overline{f}}} & \\
&\Ybs  & \rTo & \Y \\
&\dTo_{\hat{\overline{\epsilon }}_{sing}} &  & \dTo_{\overline{\epsilon}} \\
&\P^{1} &  \rTo^{t\mapsto t^{2}} & \P^{1}\\
\end{diagram}
 
\begin{remark}\label{rem: Ybarhat constructed by gluing}
Care must be taken in the choice of the conifold resolution to ensure
that the resulting space $\Ybh $ is quasi-projective (otherwise the
Gromov-Witten cannot be defined). This can be achieved by matching the
resolution obtained by constructing $\Ybh $ by gluing two copies of
$\Yh $ along $\Yh -\Sh =\hat{\epsilon }^{-1} (\C ^{*})$ via the
involution of $\Yh -\Sh $ induced by the map \eqref{eqn:
x-->x/(x^2+y^2+z^2)}. 
\end{remark}

We define the fiberwise homology groups $H_{2} (\overline{Y},\Z
)_{f}$, $H_{2} (\Ybh,\Z )_{f}$, and $H_{2} (\Ybh _{sing},\Z )_{f}$ to
be the kernels of the maps $\overline{\epsilon }_{*}$,
$\hat{\overline{\epsilon }}_{*}$, and $\left(\hat{\overline{\epsilon
}}_{sing} \right)_{*}$ respectively.

\begin{lemma}\label{lem: homology isos}
The inclusions of $S_{0}$ as the fiber over 0 and over $\infty $ in
$Y$, $\overline{Y}$, $\Ybh _{sing}$, and $\Yh _{sing}$ induce
isomorphisms in homology (taken with $\Z $ coefficients):
\[
H_{2} (S_{0})\cong H_2(Y)\cong H_2(\Ys) \cong H_{2}
(\overline{Y})_{f}\cong H_2(\Ybs)_f.
\]
Moreover, the isomorphisms induced by the inclusions of the fiber over
0 and the fiber over $\infty $ are the same isomorphism.  These groups
all have a canonical basis given by $\{C_{\rho } \}$, the components
of $\pi ^{-1} (0)$, which are labelled by $\rho \in \Irr ^{*} (G)$.

The inclusions of $\Sh $ as the fiber over 0 and over $\infty $ in
$\Ybh $ and $\Yh $ induce isomorphisms on homology (taken with $\Z $
coefficients):
\[
H_{2} (\Sh ) \cong H_2(\Yh) \cong  H_2(\Ybh)_f.
\]
Moreover, the isomorphisms induced by the inclusions of the fiber over
0 and the fiber over $\infty $ are the same isomorphism.  These groups
all have a canonical basis given by $\{\hat{C}_{\hat{\rho }} \}$, the
components of $\hat{\pi }^{-1} (0)$ which are labelled by and
$\hat{\rho }\in \Irr ^{*} (\hat{G})$.
\end{lemma}

\textsc{Proof:} Using a triangulation of the pair $(Y,\pi ^{-1} (0))$
\cite{Lojasiewicz, Hironaka} and the $\C ^{*}$ action, one can easily
construct a retract of $Y$ onto $\pi ^{-1} (0)$. It is well known that
$S_{0}$ has a deformation retract onto $\pi ^{-1} (0)$
\cite{Dimca-singularities}. Thus the inclusions
\[
\pi ^{-1} (0) \subset S_{0}\subset Y
\]
induce isomorphisms in homology. Similar arguments apply for $\pi
^{-1} (0) \subset S_{0}\subset \Yh _{sing}$ and $\hat{\pi }^{-1}
(0)\subset \Sh \subset \Yh $.

To show that $H_{2} (S_{0})\cong H_{2} (\overline{Y})_{f} $, we will
use Mayer-Vietoris for the covering $\overline{Y} = Y\cup Y$, $Y\cap
Y=Y-S_{0}$ and the compatible covering of $\P ^{1}$ to obtain the
following diagram
\[
\begin{diagram}
H_{2} (Y-S_{0})&\rTo& H_{2} (Y)\oplus H_{2} (Y)&\rTo &H_{2} (\overline{Y})&\rTo & H_{1} (Y-S_{0})\\
&& &&\dTo _{\overline{\epsilon }_{*}}&&\dTo _{\epsilon _{*}}\\
&&&&H_{2} (\P ^{1})&\rTo^{\cong} & H_{1} (\C ^{*})
\end{diagram}
\]
The Leray spectral sequence for the fibration 
\[
\epsilon :Y-S_{0}\to \C ^{*}
\]
degenerates at the $E^{2}_{pq}  $ term and gives rise to isomorphisms:
\begin{align*}
\epsilon _{*}:&H_{1} (Y-S_{0})\to H_{1} (\C ^{*})\\
i_{*}:&H_{2} (\epsilon ^{-1} (1))\to H_{2} (Y-S_{0}).
\end{align*}
Using the above isomorphisms, the commutativity of the above diagram,
and exactness, we get
\begin{align*}
H_{2} (\overline{Y})_{f}& = \ker \left(H_{2} (\overline{Y})\to H_{2}
(\P ^{1}) \right)\\
&=\ker (H_{2} \left(\overline{Y})\to H_{1} (Y-S_{0}) \right)\\
&=\coker \left(H_{2} (Y-S_{0})\to H_{2} (Y)\oplus H_{2} (Y) \right)\\
&=\coker \left(H_{2} (\epsilon ^{-1} (1))\to H_{2} (Y)\oplus H_{2} (Y) \right).
\end{align*}
By construction, the involution on $Y-S_{0}$ used to glue $Y$ to $Y$
is the identity on $\epsilon ^{-1} (1)$ so the map of $H_{2} (\epsilon
^{-1} (1))$ to $H_{2} (Y)\oplus H_{2} (Y)$ lies in the
diagonal. Consequently, the two maps $H_{2} (S_{0})\to H_{2}
(\overline{Y})_{f}$ induced by the inclusions of the fibers over 0 and
$\infty $ are the same, as asserted by the Lemma. The fact that they
are isomorphisms follows from the surjectivity of $H_{2} (\epsilon
^{-1} (1))\to H_{2} (Y)$.

The arguments for $H_{2} (S_{0})\cong H_{2} (\Ybh _{sing})_{f}$ and $H_{2} (\Sh )\cong H_{2} (\Ybh )_{f}$  are similar.\qed

\subsection{The Gromov-Witten partition function of
$\hat{Y}$}\label{subsec: GW theory of Yhat}~

\begin{prop} \label{prop: formula for Z of Yhat}
Let 
\[
\hat{c}:R^{+}\to H_{2} (\hat{S},\Z )\cong H_{2}(\hat{Y},\Z )
\]
be defined as in \S\ref{sec: intro}. Then the reduced {\GW} partition
function of $\Yh$ is given by
\[
Z^{\Yh}(\qh,\lambda)=\prod_{\alpha \in R^+}\prod _{m=1}^{\infty
}\left(1-\qh^{\,\, \hat{c} (\alpha )}\left(-e^{i\lambda} \right)^{m}
\right)^{m}.
\]
\end{prop}
\textsc{Proof:} The proof is based on the deformation invariance of
the {\GW} invariants. See \cite{BKL} and \cite{Bryan-Gholampour2} for
a similar argument.

Since $\Yh$ is smooth, $\widehat{\epsilon}$ above defines a flat
family of $\C^*$-equivariant deformations of $\Sh$. Thus it induces a
classifying morphism $$\mu: \C \rightarrow \Def(\Sh),$$ where
$\Def(\Sh)$ is the versal space of $\C^*$-equivariant deformations of
$\Sh$. $\Def(\Sh)$ is naturally identified with the complexified root
space associated to the root system $R$ (see \cite{Ka-Mo}). There is a
bijection between the set of positive roots $R^+$ and the irreducible
components of the discriminant locus in $\Def(\Sh)$. Each irreducible
component is a hyperplane in $\Def(\Sh)$ perpendicular to the
corresponding positive root.  A generic point of $\Def(\Sh)$
corresponds to an affine surface (having no compact curves), whereas a
generic point in the component of the discriminant locus perpendicular
to $\alpha \in R^+$, corresponds to a smooth surface with only one
smooth rational curve with self-intersection $-2$ representing the
homology class $\hat{c} (\alpha )$ (see \cite[Theorem~1]{Ka-Mo} or
\cite[Proposition~2.2]{BKL}).

By construction, $\mu $ is a non-constant $\C ^{*}$-equivariant map
from $\C $ to the vector space $\Def (\hat{S})$.

Since $\C^*$ acts on both the source and target with weight one (see
Remark \ref{rem:C* on epshat} and \cite[Theorem~1]{Ka-Mo}), $\mu (t)$
must be homogeneous of degree 1. Namely $\mu $ is a linear map, and
hence its image is a line passing through the origin. 

Note that because $\epsilon :Y\to \C $ has compact curves in each
fiber, the family $\hat{\epsilon }:\Yh \to \C $ has compact curves in
each fiber. Thus the image of $\mu $ is contained in the discriminant
locus.

Let 
\[
\mu _{1}:\C \to \Def (\hat{S})
\]
be a linear map whose image is a generic line passing through the
origin. Let $\hat{Y}_{1}$ be the total space of the family induced by
$\mu _{1}$. Clearly, $\hat{Y}_{1}$ is obtained from $\hat{Y}$ by a $\C
^{*}$-equivariant deformation. Since the equivariant {\GW} invariants
are invariant under equivariant deformation we have
\[
Z^{\Yh_1}(\qh,\lambda)=Z^{\Yh}(\qh,\lambda).
\]

In contrast to $\mu $, the image of $\mu_1$ meets the discriminant
locus only at the origin, hence all the curve classes of $\Yh_1$ are
represented only by the curves supported on the central fiber of the
family $\Yh_1 \rightarrow \C$. Thus, for any $\beta \neq 0$, $N_\beta
^g(\Yh_1)$ (and hence $Z^{\Yh_1}(\qh,\lambda)$) has a well defined
non-equivariant limit.  Consequently, we can compute the invariants of
$\hat{Y}_{1}$ (and hence of $\hat{Y}$) by making a non-equivariant
deformation of $\Yh_1$ to $\Yh_2$. Let
\[
\mu_2:\C \rightarrow \Def(\Sh),
\]
be a generic affine linear map whose image is an affine line meeting
the discriminant locus transversely. Let $\hat{Y}_{2}$ be the total
space of the family induced by $\mu _{2}$.  Since the {\GW} invariants
are invariant under deformations we have
\[
Z^{\Yh_2}(\qh,\lambda)=Z^{\Yh_1}(\qh,\lambda).
\]
By the description of the discriminant locus given above it is evident
that the image of $\mu_2$ meets transversally each component of the
discriminant locus once at a generic point. This implies that the
threefold $\Yh_2$ contains a unique isolated smooth rational curve
with normal bundle isomorphic to
\[
\mo(-1)\oplus \mo(-1)
\]
in the curve class corresponding to each element of $R^+$ (see
\cite[Lemma~2]{Bryan-Gholampour2}).  

The contribution to the Gromov-Witten partition function of an
isolated $(-1,-1)$ curve in a class $\beta $ is well known (e.g. the
proof of Theorem~3.1 in \cite{Behrend-Bryan}). It is given by a factor of
\[
\prod _{m=1}^{\infty }\left(1-q^{\beta }\left(-e^{i\lambda } \right)^m
\right)^{m}
\]
and thus
\[
Z^{\Yh_2}(\qh,\lambda)=\prod_{\alpha \in R^+}\prod _{m=1}^{\infty
}\left(1-\qh^{\,\,\hat{c} (\alpha )}\left(-e^{i\lambda} \right)^{m}
\right)^{m}.
\]
The proposition follows. \qed

\begin{rem}
The image of $\mu$ in $\Def (\hat{S}) $ defines a line in complex root
space and it is natural to ask which specific line it is. The
line is not quite unique---it depends on the choice of the small
resolution $\hat{Y}\to \Ys  $. We have found that there
is a choice of resolution such that the line is generated by the
vector
\[
\sum_{\hat{\rho } \in \irr(\Gh) \backslash \irr(G)} \dim (\hat{\rho} )
e_{\hat{\rho }}.
\]
The lines corresponding to other choices of the resolution are
obtained from the above line by the action of $W_{0}\subset W$, the
subgroup of the Weyl group generated by reflections about
$e_{\hat{\rho } } $, where ${\hat{\rho} }\in \Irr (\hat{G})\backslash
\Irr (G)$. The specific form of the generating vector given above was
not needed in the proof of Proposition~\ref{prop: formula for Z of
Yhat}, but it does provide a concrete description of $\hat{Y}$ and
hence $Y$ via the explicit description of the universal family over
$\Def (\hat{S}) $ given in \cite{Ka-Mo}.
\end{rem}

\subsection{Relating the Gromov-Witten theories of $Y$, $\hat{Y}$, $\overline{Y}$, and $\hat{\overline{Y}}$.}\label{subsec: relating GW of Y, Yhat, Ybar and Ybarhat}~

\begin{lemma}\label{lem: relating Z of families over C and P1}
The partition functions of $Y$ and $\hat{Y}$ are related to
the fiber class partition functions of $\overline{Y}$ and
$\hat{\overline{Y}}$ by the formulas
\begin{align*}
Z^{Y} (q,\lambda )^{2}&=Z_{f}^{\overline{Y}} (q,\lambda ),\\
Z^{\hat{Y}} (\hat{q},\lambda )^{2}&=Z_{f}^{\hat{\overline{Y}}}
(\hat{q},\lambda ).
\end{align*}
\end{lemma}

\textsc{Proof:} The $\C^*$-fixed part of the moduli space of stable
maps to $\Y$ representing fiber classes has two isomorphic components,
corresponding to maps with images in either of
$\overline{\epsilon}^{\;-1}(0)$ or
$\overline{\epsilon}^{\;-1}(\infty)$. Each component is isomorphic to
the $\C^*$-fixed part of the moduli space of stable maps to $Y$.
Since the virtual dimension of the moduli space of stable maps to $\Y$
representing $\beta \in H_2(\Y,\Z)_f$ is zero, the corresponding {\GW}
invariant $N_\beta^g(\Y)$ is independent of the weight of the $\C^*$
action. We using the isomorphisms in Lemma~\ref{lem: homology isos} to
identify the second homology groups of $\overline{\epsilon }^{-1} (0)$
and $\overline{\epsilon }^{-1} (\infty )$ with $H_{2} (\overline{Y},\Z
)_{f}$. Then by localization, we have $N_\beta^g(\Y)=2N_\beta^g(Y)$.
The same argument yields
$N_\beta^g(\hat{\overline{Y}})=2N_\beta^g(\hat{Y})$. The lemma
follows.  \qed

Finally we relate the {\GW} theory of $\hat{\overline{Y}}$ to the
{\GW} theory of $\overline{Y}$. This is given in the following
proposition:

\begin{prop} \label{prop:relation to conifold resolution} The reduced
{\GW} partition function of $\hat{\overline{Y}}$ and $\overline{Y}$
satisfy the relation
\[
\frac{Z_{f}^{\hat{\overline{Y}}}(\qh,\lambda)}{Z_{f}^{\hat{\overline{Y}}}(\qh,\lambda)|_{q=0}}
 =Z_{f}^{\overline{Y}}(q,\lambda)^2 
\]
after the specialization
\begin{equation} \label{equ:specialization to q}
\hat{q}_{\hat{\rho}} = \begin{cases}
q_{\rho } &\text{if $\,\hat{\rho }\,$ pulls back from $\rho \in \Irr ^{*} (G) $}\\
1&\text{if $\,\hat{\rho }\,$ is not the pullback of a $G$ representation.}
\end{cases}
\end{equation}
\end{prop}
\textsc{Proof:} Since $\overline{Y}$ is a family over a compact base,
and each fiber of $\overline{\epsilon}$ is a non-compact surface but
with a finite number of compact curves, the moduli space of stable
maps to $\Y$ representing a fiber class is compact. Thus,
$Z^{\Y}_f(q,\lambda)$ has a non-equivariant limit. All the elements of
$H_2(\Ybh,\Z)_f$ are Calabi-Yau classes, and by the same argument, the
corresponding {\GW} partition function, $Z^{\Ybh}_f(\qh,\lambda)$, has
a non equivariant limit.

Let 
\[
\widehat{\overline{\epsilon}}_{gen}:\Ybhg
\rightarrow \P^{1}
\]
be the double cover of the family $\overline{\epsilon}$ branched over
two generic fibers.  $\Ybhg$ is related to $\Ybh$ by
conifold transitions. Indeed, moving the branch points in $\P^{1}$ to
0 and $\infty$ defines a deformation of $\Ybhg$ to
$\Ybs$. Let $H_2(\Ybhg,\Z)_f$ be the subspace of the fiber
classes of $\Ybhg$. We have
\[
H_2(\Ybhg,\Z)_f \cong H_2(\Ybs,\Z)_f \cong
H_2(\overline{Y},\Z)_{f},
\]
where the first isomorphism is proven in \cite[page~167]{Li-Ruan} and
the second isomorphism is from Lemma~\ref{lem: homology isos}. The
{\GW} theory of $\Ybhg$ and $\Ybh$ are related by a theorem of Li and
Ruan \cite[Theorem~B]{Li-Ruan}. There is a surjective homomorphism
\[
\phi: H_2(\Ybh,\Z)_f \rightarrow H_2(\Ybhg,\Z)_f,
\]
so that the {\GW}
invariants of $\Ybhg$ and $\Ybh$ have the following
relation:
\[ 
N_\beta^g(\Ybhg)=\sum_{\beta' \in \,\phi^{-1}(\beta)}
N^g_{\beta'}(\Ybh)
\]
for any $\beta \in \, H_2(\Ybhg,\Z)_f$ with $\beta \neq 0$. Note that
by the identifications of the source and target of $\phi$ with
respectively $H_{2} (\hat{Y},\Z )$ and $H_2(Y,\Z)$ given above, $\phi$
can be identified with the homomorphism
\[
\hat{\overline{f}}_{*}:H_{2} (\hat{\overline{Y}},\Z )_{f} \rightarrow
H_2(\overline{Y},\Z)_{f}.
\]

Thus we can write
\[
\sum _{\begin{smallmatrix}\beta \in H_{2} (\overline{Y})_{f}\\
\beta \neq 0
\end{smallmatrix}} N_{\beta }^{g}
(\Ybhg ) q^{\beta } =\sum _{\begin{smallmatrix} \beta ' \in H_{2}
(\hat{\overline{Y}})_{f}\\
\hat{\overline{f}}_{*} (\beta ')\neq 0 \end{smallmatrix}} N_{\beta '}^{g} (\hat{\overline{Y}})
q^{\hat{\overline{f}}_{*} (\beta ')}.
\]
Observing that 
\[
\hat{q}^{\,\,\beta }\mapsto q^{\hat{\overline{f}}_{*} (\beta )}
\]
is the same as the specialization \eqref{equ:specialization to q}, we
obtain
\[
Z^{\Ybhg}_{f}(q,\lambda)= \frac{Z^{\Ybh}_f(\qh,\lambda)}{Z_{f}^{\Ybh }
(\hat{q},\lambda )|_{q=0}}
\]
after the specialization (\ref{equ:specialization to q}). We must
divide by the factor $Z_{f}^{\Ybh } (\hat{q},\lambda )|_{q=0}$ in
order to compensate for the lack of terms with $\hat{\overline{f}}
(\beta ')=0$ in the above sum.

It remains to be shown that 
\[
Z_{f}^{\Ybhg } (q,\lambda ) = Z_{f}^{\overline{Y}}
(q,\lambda )^{2}.
\]
The proof of this is very similar to the degeneration argument
used by Maulik and Pand\-hari\-pan\-de for the computations of the
fiber class {\GW} theory of Enriques Calabi-Yau threefold
\cite[\S1.1,\S1.4]{Maulik-Pandharipande-new} and we thank Davesh
Maulik for calling our attention to it.

Since $\Ybhg $ is given by a double cover of
$\overline{Y}$ branched over two smooth fibers, we can obtain a
degeneration of $\Ybhg $ by allowing the two
branched fibers to come together. This forms a \emph{good
degeneration} in the sense that the total space of the degeneration is
smooth and that the limiting variety is
\[
\Y\cup_F\Y,
\]
the union of two copies of $\Y$ in a normal crossing along a smooth
fiber $F$. This is the setting in which we can apply relative
Gromov-Witten theory to compute the invariants of
$\Ybhg $ in terms of the relative invariants of
the pair $(\Y,F)$.  Our situation is particularly amenable to the
computation of Gromov-Witten invariants by degeneration for two
reasons. One is that we only consider fiber curve classes and they
intersect trivially with $F$. The other is that the degeneration locus
$F$ is a (non-compact) $K3$ surface and consequently, it has no
non-trivial Gromov-Witten invariants.

For any fiber class $\beta \in
H_2(\Ybhg,\Z)_f$, and any non negative integer $g$ the
degeneration formula for {\GW} invariants implies that
\[
N_\beta^g(\Ybhg)=2N_\beta^g(\Y /F),
\]
where the right hand side is the corresponding {\GW} invariant of $\Y$
relative to $F$. However, since $F$ is a (non-compact) $K3$ surface
(see the proof of Lemma 2 in \cite{Maulik-Pandharipande-new})
\[
N_\beta^g(\Y/F)=N_\beta^g(\Y).
\]
The last two equalities immediately imply that
\[
Z^{\Ybhg}_f(q,\lambda)=Z^{\Y}_f(q,\lambda)^2.
\]
The proof of Proposition is now complete.\qed

We have now completed the proof of equations \eqref{eqn:
ZY2=ZYbar}--\eqref{eqn: ZYhat} and so the proof of
Theorem~\ref{thm:main theorem} is complete.\qed

\section{{\CRC}} \label{sec:CRC}

\subsection{Overview}

A well known principle in physics asserts that string theory on a
Calabi-Yau orbifold $X$ is equivalent to string theory on any crepant
resolution $Y\to X$. Consequently, it is expected that Gromov-Witten
theory on $Y$ should be equivalent to Gromov-Witten theory on $X$ and
this is known as the crepant resolution conjecture. The crepant
resolution conjecture goes back to Ruan \cite{Ruan-crepant} and has
recently undergone successive refinements
\cite{Bryan-Graber,CCIT,Coates-Ruan}. In the case where the orbifold
satisfies the \emph{hard Lefschetz condition}
\cite[Defn~1.1]{Bryan-Graber}, the crepant resolution conjecture says
that the Gromov-Witten potentials of $X$ and $Y$ are equal after a
change of variables. In lemma~\ref{lem: hard lefschetz => G in SO(3)
or SU(2)} we prove that the orbifold $[\C ^{3}/G]$ satisfies the hard
Lefschetz condition if and only if $G$ is a finite subgroup of $SU
(2)\subset SU (3)$ or $SO (3)\subset SU (3)$.

In this section we use the crepant resolution conjecture to derive
Conjecture~\ref{conj:orbifold prediction} and we prove
Proposition~\ref{prop: Fx conjecture implies CRC}. We will need to
work with the full genus zero Gromov-Witten potential of $Y$, which is
a generating function for the invariants of $Y$ with an arbitrary
number of insertions. Because $Y$ is Calabi-Yau, the $n$-point
Gromov-Witten invariants are easily obtained from the $0$-point
invariants along with classical intersections. The main work of this
section is computing the triple intersections on $Y$ and manipulating
the generating functions.

\subsection{Triple intersections on $Y$ and the genus zero potential
function.} Let $\{\gamma _\rho\}_{\rho \in \irr(G)}$ be the
basis for $H^2(Y,\Q)$ dual to $\{C_{\rho } \}$, namely that satisfies
\[
 \int_{C_\rho}\gamma _{\rho'}=\delta_{\rho\rho'}.
\]
Let $L_\rho$ be the line bundle on $Y$ defined by $c_1(L_\rho)=\gamma
_\rho$. We choose a lift of the $\C^*$-action of $Y$ to $L_\rho$ such
that
\[
\int_Y c_1(L_\rho)=0
\]
where the integral is evaluated via localization.  We also denote by
$\gamma_\rho$ the equivariant Chern class of $L_\rho$ with the chosen
lift. We take $\gamma_0$, corresponding to the trivial 1-dimensional
representation of $G$, to be the identity in equivariant
cohomology. We have that $\{\gamma_\rho\}_{\rho \in \Irr(G)}$ is a set
of generators for $H^*_{\C^*}(Y,\Q)$, the equivariant cohomology ring
of $Y$. Let $y=\{y_\rho\}_{\rho \in \Irr(G)}$ be a set of variables
corresponding to this cohomology basis.

For any given vector $m=(m_\rho)_{\rho \in \; \Irr(G)}$ of
nonnegative integers, we use the following notation
$$\gamma^m= \prod_{\rho \in\; \Irr(G)} \gamma_{\rho}^{m_{\rho}} \quad \text{and} \quad
\frac{y^m}{m!}=\prod_{\rho \in\; \Irr(G)}
\frac{y_{\rho}^{m_{\rho}}}{m_{\rho}!}.$$ Let $\beta \in H_2(Y,\Z)$,
$g \in \mathbb{Z}_{+}$, and suppose that
$$\displaystyle |m|=\sum_{\rho \in \Irr(G)} m_\rho.$$ The genus $g$,
degree $\beta$, $|m|$-point equivariant {\GW} invariant of $Y$
corresponding to the vector $m$ is denoted by $\langle \gamma^m
\rangle_{g,\;\beta}^{Y}.$ In general, these invariants take values in
$\Q(t)$ (see \cite{Br-Pa-local-curves}). Note that we previously
denoted the 0-point invariants, $\langle\, \rangle_{g,\;\beta}^{Y}$,
by $N^g_\beta(Y)$.
\begin{defn}
Using the notation above, we write the genus zero potential function
for {\GW} invariants of $Y$ as
\[
F^Y(q,y)=\sum_{\beta \in H_2(Y,Z)} \sum_{m \in \Z_{+}^r} \langle
\gamma^m \rangle_{0,\;\beta}^{Y} \, q^\beta\, \frac{y^m}{m!}
\]
where $r=|\Irr(G)|$. We
call the $\beta=0$ part of $F^Y$ the classical part, and the $\beta
\neq 0$ part of $F^Y$ the quantum part. We denote the former by
$F^Y_{cl}$, and the latter by $F^Y_{qu}$.
\end{defn}
$F^Y_{cl}$ involves only the classical triple equivariant
intersections on $Y$. In other words, $F^Y_{cl}$ is the cubic
polynomial in $y$ given by $$F^Y_{cl}(y)=
\sum_{\rho,\;\rho',\;\rho'' \in \; \Irr(G)} \int_Y \gamma_\rho \;
\gamma_{\rho'} \; \gamma_{\rho''} \, \,
\frac{y_\rho\;y_{\rho'}\;y_{\rho''}}{3!}.$$

If $\beta \neq 0$ then using the point and divisor axioms, one can
write all the {\GW} invariants in terms of 0-point invariants
$\langle \, \rangle_{g,\;\beta}^{Y}.$

\begin{rem} \label{rem:y=0 in quantum part}
One can see that for any $\rho \in\; \irr(G)$, the variables
$y_\rho$ and $q_\rho$ always appear in $F^Y_{qu}$ as the
product $q_\rho e^{y_\rho}$. Thus if we set $y=0$ in $F^Y$ no
information will be lost from its quantum part.
\end{rem}

We determine the classical and the quantum parts of $F^Y$ in
separate propositions.

\begin{prop} \label{prop:classical part}
The classical triple equivariant intersection numbers on $Y$ are
given by
\begin{align}
\int_Y \gamma_0 &=\frac{1}{t^3|G|}, \label{equ:0-point} \\
\int_Y \gamma_\rho &=0,\label{equ:1-point} \\
\int_Y \gamma_\rho \; \gamma_{\rho'}&=\frac{-1}{2th}\;
\sum_{\alpha \in \; R^+} \alpha^\rho\;\alpha^{\rho'},\label{equ:2-point}\\
\int_Y \gamma_\rho \;
\gamma_{\rho'}\;\gamma_{\rho''}&=\frac{1}{4}\;\sum_{\alpha \in \;
R^+}\alpha^\rho\;\alpha^{\rho'}\;\alpha^{\rho''},
\label{equ:3-point}
\end{align}
where, $t$ is the equivariant parameter, $h$ is the Coxeter number of
the root system $R$, and $\alpha ^{\rho }$ is the coefficient of the
simple root corresponding to $\rho $ in the positive root $\alpha $.
\end{prop}
\textsc{Proof:} Recall the map $\hat{f}:\hat{Y}\to Y$ constructed in
\S\ref{subsec: Geometry of Y and the fibrations}.  For any $\rho \in
\irr(G)$, the $\C^*$-equivariant line bundle $L_\rho$ on $Y$ defined
above, pulls back via $\hat{f}$ to a $\C^*$-equivariant line bundle on
$\Yh$ denoted by $\widehat{L}_{\hat{\rho }}$. By the push-pull formula
for $\hat{f}$, we have
\begin{align*}
\int_{\widehat{C}_{\hat{\rho} '}} c_1(\hat{L}_{\hat{\rho}})&=\delta_{\hat{\rho}\hat{\rho}'}\\
 \int_{\Sh} c_1(\hat{L}_{\hat{\rho}})&=0
\end{align*}
for any $\hat{\rho},\, \hat{\rho}' \in \irr(G)\subset \irr
(\hat{G})$. We extend $\{\widehat{L}_{\hat{\rho}}\}_{\hat{\rho} \in
\irr(G)}$ to the set of $\C^*$-equivariant line bundles
$\{\widehat{L}_{\hat{\rho}}\}_{\hat{\rho} \in \irr(\Gh)}$ satisfying the
two relations above for any $\hat{\rho},\, \hat{\rho}' \in
\irr(\hat{G})$. Define
$\widehat{\gamma}_{\hat{\rho}}=c_1(\widehat{L}_{\hat{\rho}})$, and let
$\widehat{\gamma}_0$ be the class of identity in $H^*_{\C^*}(\Sh,\Q)$,
the equivariant cohomology of $\Sh$. By construction,
$\{\widehat{\gamma}_{\hat{\rho}}\}_{\hat{\rho} \in \Irr(\Gh)}$ is a
basis for $H^*_{\C^*}(\Sh,\Q)$.

Now let $\widehat{\gamma} \in H_{\C ^{*}}^*(\Yh,\Q)$ be the pull back
of some $\gamma \in H_{\C ^{*}}^*(Y,\Q)$ via the $\C^*$-equivariant
map $\hat{f}$. We can write
\[
\int_Y \gamma = \frac{1}{2} \int_{\Yh}\widehat{\gamma}=
\frac{1}{2t}\int_{\Sh}\widehat{\gamma}.
\]
All the integrals are defined via localization. The factor 1/2 in the
first equality appears because $\hat{f}$ is generically a double cover
map, and the $1/t$ factor in the second equality because the weight of
the normal bundle of $\Sh\hookrightarrow \Yh$ is 1. Proposition
\ref{prop:classical part} follows immediately from these identities
and Lemma \ref{lem:integrals on surface}. \qed

\begin{lem} \label{lem:integrals on surface} 
Let $\{\widehat{\gamma}_{\hat{\rho}}\}_{\hat{\rho} \in \Irr(\Gh)}$ be the basis for
$H^*_{\C^*}(\Sh,\Q)$ chosen above. Then we have
\begin{align*}
\int_{\Sh} \widehat{\gamma}_0 &=\frac{4}{t^2|\Gh|},\\
\int_{\Sh} \widehat{\gamma}_{\hat{\rho}} &=0,\\
\int_{\Sh} \widehat{\gamma}_{\hat{\rho}} \; \widehat{\gamma}_{{\hat{\rho}}'}&=\frac{-1}{h}
\sum_{\alpha \in \; R^+}\alpha^{\hat{\rho}}\;\alpha^{{\hat{\rho}}'},\\
\int_{\Sh} \widehat{\gamma}_{\hat{\rho}} \;
\widehat{\gamma}_{{\hat{\rho}}'}\;\widehat{\gamma}_{{\hat{\rho}}''}&=\frac{t}{2}\sum_{\alpha
\in \; R^+}\alpha^{\hat{\rho}}\;\alpha^{{\hat{\rho}}'}\;\alpha^{{\hat{\rho}}''},
\end{align*}
\end{lem}
\textsc{Proof:} This is basically \cite[Lemma~4]{Bryan-Gholampour2}
with only the following differences 
\begin{enumerate} [1.] \item The basis chosen for $H^*_{\C^*}(\Sh,\Q)$
is obtained from the basis chosen in \cite{Bryan-Gholampour2} by the
negative inverse of the Cartan matrix of $R$. The reason for this is
that the cohomology basis chosen in \cite{Bryan-Gholampour2} is the
dual of $\{C_{\hat{\rho}}\}_{{\hat{\rho}} \in \irr(\Gh)}$ with respect to the
intersection pairing (which is the negative Cartan matrix by the
{\MC}), while in this paper it is the dual of $\{C_{\hat{\rho}}\}_{{\hat{\rho}} \in
\irr(\Gh)}$ with respect to the equivariant Poincar\'e pairing.
\item The weights of the induced $\C^*$-action on the canonical
bundle of $\Sh$ and the cohomology basis elements
$\widehat{\gamma}_{\hat{\rho}}$ are half of the corresponding weights in
\cite{Bryan-Gholampour2}.
\end{enumerate} \qed

\begin{rem}
Note that (\ref{equ:1-point}) in Proposition \ref{prop:classical
part} is the direct result of the definition of $\{\gamma_\rho\}$.
(\ref{equ:0-point}) in Proposition \ref{prop:classical part} can
also be deduced by the same argument as in the proof of Lemma 4 in
\cite{Bryan-Gholampour2}, by pushing forward the class of the identity to
$H^0_{\C^*}(\X)$.
\end{rem}

\begin{rem} \label{rem:paring}
Define the matrix $G=(g_{\rho \rho'})_{\rho,\,\rho' \in\, \irr(G)}$
by $$g_{\rho \rho'}= \int_Y \gamma_\rho \; \gamma_{\rho'}.$$ Then
one can show that $G^{-1}=(g^{\rho \rho'})$ with
\[
g^{\rho \rho'}= t\left\langle (V-3\C)\otimes \rho, \rho'  \right\rangle
\]
where $\LangRang{\cdot ,\cdot }$ is the pairing making the irreducible
representations orthonormal. The matrix $G^{-1}$ is the McKay quiver
and is the analog of the negative Cartan matrix in dimension 2.
\end{rem}

\begin{prop}
The genus zero {\GW} potential function of $Y$ is given by
\[
F^Y(q,y)= F^Y_{cl}(y)+ \sum_{\alpha \in \;
R^+}\sum^{\infty}_{d=1}\frac{1}{2d^3}\,q^{dc (\alpha )}e^{d\sum _{\rho }\alpha ^{\rho }y_{\rho }}.
\]
\end{prop}
\textsc{Proof:} This is a consequence of Theorem~\ref{thm:main
theorem} and Remark~\ref{rem:y=0 in quantum part}. Each genus zero
Gopakumar-Vafa invariant $n_{\beta }^{0}$ contributes ${n_{\beta
}^{0}}{d^{-3}}$ to the degree $d\beta $ genus zero Gromov-Witten
invariant. The proposition follows.\qed 
\smallskip

\subsection{Deriving $F^{\X }$ via the crepant resolution conjecture.}~

Let $\con^*(G)$ be the set of nontrivial conjugacy classes of $G$. For
$(g)\in \con ^{*} (G )$ and $\rho \in \irr (G)$, we define
\[ 
L_\rho^{\g}=\frac{|\,\text{C}(g)|}{|G|} \sqrt{3-\chi_V(g)}\;\chi_\rho ^{g},
\]
where $\text{C}(g)$ is the centralizer of the group element $g$, and
$V$ is the 3 dimensional representation of $G$ induced by the
embedding $G\subset SO (3)$. It was conjectured in
\cite[Conjecture~3.1]{Bryan-Graber} that $F^\X(x)$ is obtained from
$F^Y (y,q)$ after replacing

\begin{align}\label{eqn: change of vars for CRC}
\nonumber y_0 &=x_{(e)},  \\ 
y_\rho &=\sum_{(g) \in\;
\con^*(G)}i\,L_\rho^{\g} \, x_{\g},\\
\nonumber q_\rho &=\exp \left(\frac{2\pi i\dim (\rho )}{|G|} \right).
\end{align}

Recalling the definition of $\mathsf{X}_{\rho }$ given in
conjecture~\ref{conj:orbifold prediction}, we see that
\[
\mathsf{X}_{\rho } = \frac{2\pi \dim (\rho )}{|G|}+\sum _{(g)\in \con
^{*} (G)} L_{\rho }^{(g)}x_{(g)}.
\] 
Using the above equations to substitute the $y$ variables for the $x$
variables, and specializing the $q$ variables, we obtain the prediction
\begin{align*}
F^{\X } (x)= F^Y_{cl}(x)+\sum_{\alpha \in\; R^+} \frac{1}{2d^3}
\exp\left(id\sum _{\rho }\alpha ^{\rho }\mathsf{X}_{\rho } \right).
\end{align*}
We take the third partial derivatives of the above equation with
respect to the variables $x_{(k)}$, $x_{(k')}$ and $x_{(k'')}$
corresponding to three nontrivial conjugacy classes, use (4) in
Proposition \ref{prop:classical part}, and sum the geometric series to
arrive at
\begin{gather*}
F^{\X }_{x_{(k)}x_{(k')}x_{(k'')}}(x) 
=\sum_{\rho,\,\rho',\,\rho''}\bigg(\frac{1}{4}\sum_{\alpha \in \;
R^+}\alpha^{\rho}\alpha^{\rho'}\alpha^{\rho''}\bigg)i^3\,L_\rho^{(k)}\,L_{\rho'}^{(k')}\,L_{\rho''}^{(k'')}
\\ +\sum_{\alpha\in\; R} \frac{i^3}{2}\; \bigg(\sum_\rho
\alpha^{\rho}L_\rho^{(k)}\cdot\sum_{\rho'}
\alpha^{\rho'}L_{\rho'}^{(k')}\cdot\sum_{\rho''}
\alpha^{\rho''}L_{\rho''}^{(k'')}\bigg)\; \frac{e^{i\sum _{\rho
}\alpha ^{\rho }\mathsf{X}_{\rho }}}{1-e^{i\sum _{\rho }\alpha ^{\rho
}\mathsf{X}_{\rho }}}. \nonumber
\end{gather*}

Using the identity
\[
\frac{e^{i\;\theta}}{1-e^{i\;\theta}}=-\frac{1}{2}-\frac{i}{2}\tan\left(\frac{\theta
}{2} +\frac{\pi}{2}\right),
\]
we find that a delicate cancellation leads to

\begin{gather*}
F^{\X }_{x_{(k)}x_{(k')}x_{(k'')}}(x) =\\
\frac{-1}{4}\sum_{\alpha \in\; R^+} \bigg(\sum_\rho
\alpha^{\rho}L_\rho^{(k)}\cdot \sum_{\rho'}
\alpha^{\rho'}L_{\rho'}^{(k')}\cdot \sum_{\rho''}
\alpha^{\rho''}L_{\rho''}^{(k'')}\bigg)\tan\left(\frac{1}{2}\sum _{\rho }\alpha ^{\rho
}\mathsf{X}_{\rho } +\frac{\pi}{2}\right).
\end{gather*}
By integrating the above with respect to the variables
$x_{(k)}$, $x_{(k')}$ and $x_{(k'')}$, we  have derived the prediction in
Conjecture \ref{conj:orbifold prediction} for $F^\X (x)$.

\begin{rem} \label{rem:x0 part}
By the point axiom, the only nontrivial terms in $F^\X(x)$
containing $x_{(e)}$ are $$\langle\, \delta_{(e)}^3 \, \rangle
\,\frac{x_{(e)}^3}{3!} \quad \text{ and } \quad \langle
\,\delta_{(e)}\delta_{(g)}\delta_{(g^{-1})}\, \rangle\,
x_{(e)}x_{\g}x_{(g^{-1})} \text{ for } g\neq e$$ corresponding to
the classical invariants, which are easily evaluated to
$$\langle\, \delta_{(e)}^3 \, \rangle=\frac{1}{t^3|G|} \quad
\text{ and } \quad \langle
\,\delta_{(e)}\delta_{(g)}\delta_{(g^{-1})}\,
\rangle=\frac{1}{t|\text{C}(g)|}.$$ By the similar argument given in
\cite{Bryan-Gholampour2}, one can check that these terms match up with the part
of $F^Y$ with nonzero $y_0$ terms (see Remark \ref{rem:paring}).
\end{rem}

The above derivation of $F^{\X } (x)$, along with the preceding remark
shows that if Conjecture~\ref{conj:orbifold prediction} is true, then
the Crepant Resolution Conjecture holds for $Y\to \C^3/G$ using the
change of variables given by equation~\eqref{eqn: change of vars for
CRC}. This proves Proposition~\ref{prop: Fx conjecture implies CRC}.

\subsection{The Hard Lefschetz condition}\label{subsec: hard lefschetz}~

\begin{lemma}\label{lem: hard lefschetz => G in SO(3) or SU(2)}
Let $G$ be a finite subgroup of $SU (3)$, then the following
conditions are equivalent:
\begin{enumerate}
\item The orbifold $[\C ^{3}/G]$ satisfies the hard Lefschetz
condition \cite[Defn~1.1]{Bryan-Graber}.
\item Every non-trivial element of $G$ has age 1.
\item $G$ is a finite subgroup of $SO (3)\subset SU (3)$ or $SU
(2)\subset SU (3)$.
\end{enumerate}
\end{lemma}
\textsc{Proof:} Let $g\in G$ be an element of order $n>1$. The
eigenvalues of the action of $g$ on $\C ^{3}$ are $\omega ^{k_{1}}$,
$\omega ^{k_{2}}$, and $\omega ^{k_{3}}$, where $\omega
=\exp\left(\frac{2\pi i}{n} \right)$, $k_{i}\in \{0,\dotsc ,n-1 \}$,
and $k_{1}+k_{2}+k_{3}\equiv 0\mod n$. By definition, the age of the
element $g$ is given by
\[
\age (g) = \frac{1}{n} (k_{1}+k_{2}+k_{3}).
\]
By definition, the hard Lefschetz condition means that $\age (g)=\age
(g^{-1})$ for all $g$. If $k_{i}\neq 0$ for $i=1,2,3$, then 
\[
\age (g^{-1}) = \frac{1}{n} ((n-k_{1})+ (n-k_{2})+ (n-k_{3})) = 3-\age (g).
\]
Consequently, if $\age (g)=\age (g^{-1})$ then $k_{i}=0$ for some
$i$. In this case the eigenvalues must be $(1,\omega ^{k},\omega
^{n-k})$ and $\age (g)=\age (g^{-1})=1$. This proves that (1) is
equivalent to (2).

If $G$ is a subgroup of $SU (2)\subset SU (3)$ or $SO (3)\subset SU
(3)$, then every element fixes some line in $\C ^{3}$ and hence has
age one (or zero). Thus (3) implies (2), and it remains to be shown
that (2) implies (3).

Let $\chi _{V}$ be the character of $V$, the 3 dimensional
representation induced by $G\subset SU (3)$. By assumption, the action
of $g$ on $V$ has eigenvalues $\{1,\omega ^{k},\omega ^{-k} \}$. We have
\begin{align*}
\chi _{V} (g^{2})&= 1+\omega ^{2k}+\omega ^{-2k}\\
&= (\omega ^{k}+\omega ^{-k}+1) (\omega ^{k}+\omega ^{-k}-1)\\
&=\overline{\chi} _{V} (g)\left(\chi _{V} (g)-2 \right).
\end{align*}
Thus
\[
\frac{1}{|G|}\sum _{g\in G} \chi _{V} (g^{2}) = \frac{1}{|G|}\sum
_{g\in G} \overline{\chi }_{V} (g)\left(\chi _{V} (g)-2 \right).
\]
We decompose $V$ into irreducible summands and we let $r$, $c$, and
$q$ be the number of real, complex, and quaternionic summands
respectively and let $t$ be the number of trivial summands. By
\cite[Ex.~3.38]{Fulton-Harris}, the left hand side of the above
equation is equal to $r-q$ while by the orthogonality of characters,
the right hand side is equal to $r+c+q-2t$. Thus we obtain
\[
2t=c+2q
\]
which, since $V$ is of dimension 3, can only be solved by $t=0$ or
$t=1$. If $t=0$ then $c=q=0$, $r=1$, and so $V$ is real,
i.e. $G\subset SO (3)$. If $t=1$ then $G\subset SU (2)$ and then lemma
is proved.\qed

\section{The $D_{5}$ case: $G=\St $.}\label{sec: D5 case}~

To illustrate the theorems of this paper, we work out the $D_{5}$ case
explicitly.  Let $\St \subset SO(3)$ be the group of permutations of
three letters (which is the dihedral group of the triangle and
corresponds to $D_{5}$ in the classification).  We have
\[
\irr(\St)=\{V_1,V_2\}
\]
where $V_1$ and $V_2$ are respectively the alternating and the
standard representations of $\St$. The embedding $\St\subset SO (3)$
induces the 3-dimensional representation
\[
V=V_1\oplus V_2.
\]

The binary version of $\St$ is the generalized quaternion group of
order 12 which we denote by $\Sth$. The non-trivial irreducible
representations of $\Sth $ are
\[
\irr(\Sth)=\{V_1,V_2,U_1,U_2,U_3\}.
\]
$\Sth$ acts on $\C^2$ via $U_{1}$. Recall that
\[
\pih: \Sh \rightarrow \C^2/\Sth\quad \quad \pi :Y\to \C ^{3}/\St 
\]
are the minimal and $G$-Hilbert scheme resolutions. The intersection
graphs of the exceptional sets $\pih^{-1}(0) $ and $\pi^{-1}(0) $ are
given in Figure~\ref{fig:intersection graphs of the exceptional
sets}. The vertices in the graphs are corresponding to the irreducible
components of $\pih^{-1}(0) $ and $\pi^{-1}(0) $. In each graph two
vertices are connected if and only if the corresponding components
meet at one point. All the components of $\pih^{-1}(0) $ are smooth
rational curves having normal bundles isomorphic to $\mo (-2)$. The
components of $\pi^{-1}(0) $ corresponding to $V_1$ and $V_2$ are
smooth rational curves with normal bundles respectively isomorphic to
\[
\mo(-1)\oplus \mo(-1)\quad  \text{and} \quad \mo(1)\oplus
\mo(-3).
\]
\begin{center} \begin{figure}[ht]
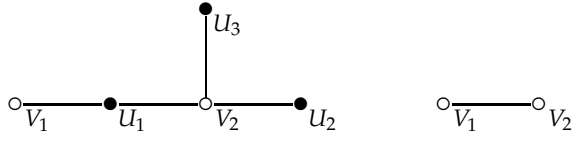
 \centering
\begin{diagram}[height=1.5em,width=1.5em,abut]
      &       &        &       &    \bullet \putlabelUUU   & &   &&&      &        &     \\
      &       &        &       &\dLine &    &   &&&      &        &     \\
\circ \putlabelV &\hLine &\bullet \putlabelU &\hLine &\circ
\putlabelVV  & \hLine & \bullet \putlabelUU    &&&
\circ  \putlabelVVV &\hLine  &\circ \putlabelVVVV \\
\end{diagram}
\caption{The intersection graphs of $\pih^{-1}(0) $ and
$\pi^{-1}(0) $ respectively.} \label{fig:intersection graphs
of the exceptional sets}
\end{figure}
\end{center}
\afterpage{} The morphism $f:\Sh \rightarrow Y$ contracts the
components corresponding to the solid vertices. Let $C_{i} \in
H_2(Y,\Z)$ be the curve class corresponding to $V_i$ for $i=1,\,2$. We
denote the homology class $d_1C_{1}+d_2C_{2}$ by $(d_1,d_2)$.

\begin{remark}\label{rem: normal bundles}
In general, the morphism $f:\hat{S}\to Y$ contracts disjoint rational
curves (see Figure~\ref{fig:ADE}). The normal
bundle of a curve $C$ in $\pi ^{-1} (0)$ is isomorphic to $\mo
(-k)\oplus \mo (k-2)$ where $k$ is the number of components collapsed
by $f$ onto $C$. This is a consequence of
\cite[Theorem~3.2]{Morrison}.
\end{remark}

The left graph in Figure~\ref{fig:intersection graphs of the
exceptional sets}, is the Dynkin diagram of the $\text{D}_5$ root
system and has vertices in correspondence with the simple roots. There
are 20 positive roots in this root system (the first three are the
binary roots): \medskip
\begin{center}
\small{
\begin{tabular}{@{}l@{}l@{}l@{}l@{}l@{}l@{}l@{}l@{}l@{}l@{}}
 \mbox{\renewcommand{\arraystretch}{.5} \begin{tabular}{l@{}l@{}l@{}l@{}}
     & &0&\\
   0&0&0&1 
   \end{tabular}}&
\mbox{\renewcommand{\arraystretch}{.5}
\begin{tabular}{l@{}l@{}l@{}l@{}}
     & &0&\\
    0&1&0&0 \;
   \end{tabular}}&
\mbox{\renewcommand{\arraystretch}{.5}
\begin{tabular}{l@{}l@{}l@{}l@{}}
     & &1&\\
    0&0&0&0 \;
   \end{tabular}}&
\mbox{\renewcommand{\arraystretch}{.5}
\begin{tabular}{l@{}l@{}l@{}l@{}}
     & &0&\\
    1&0&0&0 \;
   \end{tabular}}&
\mbox{\renewcommand{\arraystretch}{.5}
\begin{tabular}{l@{}l@{}l@{}l@{}}
     & &0&\\
    0&0&1&0 \;
   \end{tabular}}&
    \mbox{\renewcommand{\arraystretch}{.5} \begin{tabular}{l@{}l@{}l@{}l@{}}
     & &0&\\
    1&1&0&0 \;
   \end{tabular}}&
\mbox{\renewcommand{\arraystretch}{.5}
\begin{tabular}{l@{}l@{}l@{}l@{}}
     & &0&\\
    0&1&1&0 \;
   \end{tabular}}&
\mbox{\renewcommand{\arraystretch}{.5}
\begin{tabular}{l@{}l@{}l@{}l@{}}
     & &0&\\
    0&0&1&1 \;
   \end{tabular}}&
\mbox{\renewcommand{\arraystretch}{.5}
\begin{tabular}{l@{}l@{}l@{}l@{}}
     & &1&\\
    0&0&1&0 \;
   \end{tabular}}&
\mbox{\renewcommand{\arraystretch}{.5}
\begin{tabular}{l@{}l@{}l@{}l@{}}
     & &0&\\
    1&1&1&0 \;
   \end{tabular}}\\
\quad \\
\mbox{\renewcommand{\arraystretch}{.5}
\begin{tabular}{l@{}l@{}l@{}l@{}}
     & &0&\\
   0&1&1&1 
   \end{tabular}}&
\mbox{\renewcommand{\arraystretch}{.5}
\begin{tabular}{l@{}l@{}l@{}l@{}}
     & &1&\\
    0&1&1&0 
   \end{tabular}}&
\mbox{\renewcommand{\arraystretch}{.5}
\begin{tabular}{l@{}l@{}l@{}l@{}}
     & &1&\\
    0&0&1&1 
   \end{tabular}}&
\mbox{\renewcommand{\arraystretch}{.5}
\begin{tabular}{l@{}l@{}l@{}l@{}}
     & &0&\\
    1&1&1&1 
   \end{tabular}}&
\mbox{\renewcommand{\arraystretch}{.5}
\begin{tabular}{l@{}l@{}l@{}l@{}}
     & &1&\\
    1&1&1&0 
   \end{tabular}}&
    \mbox{\renewcommand{\arraystretch}{.5} \begin{tabular}{l@{}l@{}l@{}l@{}}
     & &1&\\
    0&1&1&1 
   \end{tabular}}&
\mbox{\renewcommand{\arraystretch}{.5}
\begin{tabular}{l@{}l@{}l@{}l@{}}
     & &1&\\
    0&1&2&1 
   \end{tabular}}&
\mbox{\renewcommand{\arraystretch}{.5}
\begin{tabular}{l@{}l@{}l@{}l@{}}
     & &1&\\
    1&1&1&1 
   \end{tabular}}&
\mbox{\renewcommand{\arraystretch}{.5}
\begin{tabular}{l@{}l@{}l@{}l@{}}
     & &1&\\
    1&1&2&1 
   \end{tabular}}&
\mbox{\renewcommand{\arraystretch}{.5}
\begin{tabular}{l@{}l@{}l@{}l@{}}
     & &1&\\
    1&2&2&1
   \end{tabular}}\\
\end{tabular}
}
\end{center}
\medskip where the numbers are the coefficients of the corresponding
simple root in the graph above.

By Theorem \ref{thm:main theorem}, the partition function for the
{\GW} invariants of $Y$ is given by
\[ 
Z^Y(q_{1},q_{2},\lambda)=
M(q_{1})M(q_{1}q_{2})^{2}
M(q_{2})^{4}M(q_{2}^2)^{1/2}M(q_{1}q_{2}^2),
\]
where 
\[
M (q) = \prod _{m=1}^{\infty }\left(1-q\left(-e^{i\lambda }
\right)^{m} \right)^{m}.
\]

The nonzero BPS states are
\[
n^0_{(1,0)}=1,  \quad n^0_{(1,1)}=2, \quad n^0_{(0,1)}=4, \quad n^0_{(0,2)}=1/2, \quad n^0_{(1,2)}=1.
\]
In particular, we obtain a somewhat curious genus zero multiple cover
formula for the rational curve $C_{2}$, which has normal bundle $\mo
(1)\oplus \mo (-3)$:
\[
N^{0}_{dC_{2}} = \begin{cases}
\frac{4}{d^{3}} &\text{if $d$ is odd,}\\
\frac{8}{d^{3}} &\text{if $d$ is even.}
\end{cases}
\]
Here the $\C ^{* }$ action acts trivially on $C_{2}$ and acts on the
normal bundle with weights 2 on the $\mo (-3)$ factor and 1 on the
$\mo (1)$ factor.

Conjecture~\ref{conj:orbifold prediction} can be written out
explicitly in this case. Let $x_{1}$ and $x_{2}$ be the variables
corresponding to the conjugacy class of the transpositions and the
three cycles respectively. Then the genus zero orbifold potential
function of $[\C ^{3}/\St ]$ is predicted to be
\begin{align*}
F (x_{1},x_{2}) =&\quad \,  \mathbf{h}\left( \frac{4\pi }{3}-x_{1}+\frac{x_{2}}{\sqrt{3}}\right)+4\mathbf{h}\left(\frac{\pi }{3}+\frac{x_{2}}{\sqrt{3}} \right)+2\mathbf{h}\left(x_{1} \right) \\
&+\mathbf{h}\left( \frac{4\pi }{3}+x_{1}+\frac{x_{2}}{\sqrt{3}}\right)+\frac{1}{2}\mathbf{h}\left(\frac{\pi }{3}-\frac{2x_{2}}{\sqrt{3}} \right)\\
&\quad \\
=&\quad \, \frac{1}{2}x_{1}^{2}x_{2}+\frac{1}{18}x_{2}^{3} -\frac{5}{48}x_{1}^{4}-\frac{1}{6}x_{1}^{2}x_{2}^{2}-\frac{1}{36}x_{2}^{4} +\\
&+\frac{1}{12}x_{1}^{4}x_{2}+\frac{1}{18}x_{1}^{2}x_{2}^{3}+\frac{1}{324}x_{2}^{5}+\dotsb 
\end{align*}

Since (most of the components of) the moduli space of genus zero maps
to $[\C^3/\St ]$ can be identified with the space of degree three
admissible covers, the orbifold invariants predicted above can also be
viewed as integrals of certain Hodge classes over the moduli space of
degree three admissible covers. For example, the special cases where
there are exactly two simple ramifications and the rest are double
ramification was studied in \cite{Bryan-Graber-Pandharipande}. The
series $B (u)$ which is defined and computed in
\cite[Appendix~A]{Bryan-Graber-Pandharipande} is given by $F_{112}
(0,-u)$ and indeed, the formula
\[
B (u) = \frac{1}{\sqrt{3}}\tan \left(\frac{u}{\sqrt{12}}+\frac{\pi }{3} \right)
\]
given in \cite{Bryan-Graber-Pandharipande} can be easily recovered
from the above prediction for $F (x_{1},x_{2})$.

\subsection{Acknowledgements.}  The authors thank the referees for
valuable comments on both technical and expository aspects of this
paper. The authors also warmly thank Davesh Maulik, Miles Reid, Rahul
Pandharipande, Sophie Terouanne, Michael Thaddeus, and Hsian-Hua Tseng
for helpful conversations. The authors also acknowledge support from
NSERC, MSRI, the Killam Trust, and the Miller Institute.

\bibliography{mainbiblio}

\begin{thebibliography}{10}

\bibitem{Ab-Gr-Vi}
Dan Abramovich, Tom Graber, and Angelo Vistoli.
\newblock Algebraic orbifold quantum products.
\newblock In {\em Orbifolds in mathematics and physics (Madison, WI, 2001)},
  volume 310 of {\em Contemp. Math.}, pages 1--24. Amer. Math. Soc.,
  Providence, RI, 2002.

\bibitem{AKMV}
Mina Aganagic, Albrecht Klemm, Marcos Mari{\~n}o, and Cumrun Vafa.
\newblock The topological vertex.
\newblock {\em Comm. Math. Phys.}, 254(2):425--478, 2005.
\newblock arXiv:hep-th/0305132.

\bibitem{Behrend-Bryan}
Kai Behrend and Jim Bryan.
\newblock {Super-rigid Donaldson-Thomas invariants}.
\newblock {\em Mathematical Research Letters}, 14(4):559--571, 2007.
\newblock arXiv version: math.AG/0601203.

\bibitem{Boissiere-Sarti}
Samuel Boissiere and Alessandra Sarti.
\newblock {Contraction of excess fibres between the McKay correspondences in
  dimensions two and three}.
\newblock {\em Ann. Inst. Fourier}, 57(6):1839--1861, 2007.
\newblock arXiv version: math.AG/0504360.

\bibitem{BKR}
Tom Bridgeland, Alastair King, and Miles Reid.
\newblock The {M}c{K}ay correspondence as an equivalence of derived categories.
\newblock {\em J. Amer. Math. Soc.}, 14(3):535--554 (electronic), 2001.

\bibitem{Bryan-Gholampour2}
Jim Bryan and Amin Gholampour.
\newblock {Root systems and the quantum cohomology of ADE resolutions}.
\newblock {\em Algebra and Number Theory}, 2(4):369--390, 2008.
\newblock arXiv:0707.1337.

\bibitem{Bryan-Gholampour1}
Jim Bryan and Amin Gholampour.
\newblock {Hurwitz-Hodge integrals, the $E_{6}$ and $D_{4}$ root systems, and
  the Crepant Resolution Conjecture}.
\newblock {\em Advances in Mathematics}, 221(4):1047--1068, 2009.
\newblock arXiv:0708.4244.

\bibitem{Bryan-Graber}
Jim Bryan and Tom Graber.
\newblock The crepant resolution conjecture.
\newblock In {\em Algebraic Geometry---Seattle 2005}, volume~80 of {\em Proc.
  Sympos. Pure Math.}, pages 23--42. Amer. Math. Soc., Providence, RI, 2009.
\newblock arXiv: math.AG/0610129.

\bibitem{Bryan-Graber-Pandharipande}
Jim Bryan, Tom Graber, and Rahul Pandharipande.
\newblock The orbifold quantum cohomology of {$\mathbf{C^2/Z_3}$ and Hurwitz
  Hodge integrals}.
\newblock {\em Journal of Algebraic Geometry}, 17:1--28, 2008.
\newblock arXiv version:math.AG/0510335.

\bibitem{BKL}
Jim Bryan, Sheldon Katz, and Naichung~Conan Leung.
\newblock Multiple covers and the integrality conjecture for rational curves in
  {C}alabi-{Y}au threefolds.
\newblock {\em J. Algebraic Geom.}, 10(3):549--568, 2001.
\newblock Preprint version: math.AG/9911056.

\bibitem{Br-Pa}
Jim Bryan and Rahul Pandharipande.
\newblock B{P}{S} states of curves in {C}alabi-{Y}au 3-folds.
\newblock {\em Geom. Topol.}, 5:287--318 (electronic), 2001.
\newblock arXiv: math.AG/0009025.

\bibitem{Br-Pa-local-curves}
Jim Bryan and Rahul Pandharipande.
\newblock {The local Gromov-Witten theory of curves}.
\newblock {\em Journal of the American Mathematical Society}, 21:101--136,
  2008.
\newblock arXiv:math.AG/0411037.

\bibitem{Cavalieri-Hurwitz-Hodge}
Renzo Cavalieri.
\newblock {Generating Functions for Hurwitz-Hodge Integrals}.
\newblock math.AG/math/0608590.

\bibitem{Chen-Ruan}
Weimin Chen and Yongbin Ruan.
\newblock Orbifold {G}romov-{W}itten theory.
\newblock In {\em Orbifolds in mathematics and physics (Madison, WI, 2001)},
  volume 310 of {\em Contemp. Math.}, pages 25--85. Amer. Math. Soc.,
  Providence, RI, 2002.

\bibitem{Chiodo}
Alessandro Chiodo.
\newblock {Towards an enumerative geometry of the moduli space of twisted
  curves and $r$-th roots}.
\newblock arXiv:math/0607324.

\bibitem{CCIT-CRC}
Tom Coates, Alessio Corti, Hiroshi Iritani, and Hsian-Hua Tseng.
\newblock {The Crepant Resolution Conjecture for Type A Surface Singularities}.
\newblock arXiv:0704.2034v1 [math.AG].

\bibitem{CCIT}
Tom Coates, Alessio Corti, Hiroshi Iritani, and Hsian-Hua Tseng.
\newblock {Wall-Crossings in Toric Gromov-Witten Theory I: Crepant Examples}.
\newblock arXiv:math.AG/0611550.

\bibitem{Coates-Ruan}
Tom Coates and Yongbin Ruan.
\newblock {Quantum Cohomology and Crepant Resolutions: A Conjecture}.
\newblock arXiv:0710.5901.

\bibitem{coxkatz}
David~A. Cox and Sheldon Katz.
\newblock {\em Mirror symmetry and algebraic geometry}.
\newblock American Mathematical Society, Providence, RI, 1999.

\bibitem{Dimca-singularities}
Alexandru Dimca.
\newblock {\em Singularities and topology of hypersurfaces}.
\newblock Universitext. Springer-Verlag, New York, 1992.

\bibitem{Fulton-Harris}
William Fulton and Joe Harris.
\newblock {\em Representation theory}.
\newblock Springer-Verlag, New York, 1991.
\newblock A first course, Readings in Mathematics.

\bibitem{Go-Va}
Rajesh Gopakumar and Cumrun Vafa.
\newblock M-theory and topological strings--{II}, 1998.
\newblock Preprint, hep-th/9812127.

\bibitem{Hironaka}
Heisuke Hironaka.
\newblock Triangulations of algebraic sets.
\newblock In {\em Algebraic geometry ({P}roc. {S}ympos. {P}ure {M}ath., {V}ol.
  29, {H}umboldt {S}tate {U}niv., {A}rcata, {C}alif., 1974)}, pages 165--185.
  Amer. Math. Soc., Providence, R.I., 1975.

\bibitem{mirrorbook}
Kentaro Hori, Sheldon Katz, Albrecht Klemm, Rahul Pandharipande, Richard
  Thomas, Cumrun Vafa, Ravi Vakil, and Eric Zaslow.
\newblock {\em Mirror symmetry}, volume~1 of {\em Clay Mathematics Monographs}.
\newblock American Mathematical Society, Providence, RI, 2003.
\newblock With a preface by Vafa.

\bibitem{Ka-Mo}
Sheldon Katz and David~R. Morrison.
\newblock Gorenstein threefold singularities with small resolutions via
  invariant theory for {W}eyl groups.
\newblock {\em J. Algebraic Geom.}, 1(3):449--530, 1992.

\bibitem{Li-Ruan}
An-Min Li and Yongbin Ruan.
\newblock Symplectic surgery and {G}romov-{W}itten invariants of {C}alabi-{Y}au
  3-folds.
\newblock {\em Invent. Math.}, 145(1):151--218, 2001.

\bibitem{Li-Liu-Liu-Zhou}
Jun Li, Chiu-Chu~Melissa Liu, Kefeng Liu, and Jian Zhou.
\newblock {A Mathematical Theory of the Topological Vertex}.
\newblock arXiv:math/0408426.

\bibitem{Lojasiewicz}
S.~Lojasiewicz.
\newblock Triangulation of semi-analytic sets.
\newblock {\em Ann. Scuola Norm. Sup. Pisa (3)}, 18:449--474, 1964.

\bibitem{MNOP1}
D.~Maulik, N.~Nekrasov, A.~Okounkov, and R.~Pandharipande.
\newblock Gromov-{W}itten theory and {D}onaldson-{T}homas theory. {I}.
\newblock {\em Compos. Math.}, 142(5):1263--1285, 2006.
\newblock arXiv:math.AG/0312059.

\bibitem{Maulik-Pandharipande-new}
Davesh Maulik and Rahul Pandharipande.
\newblock {New calculations in Gromov-Witten theory}.
\newblock arXiv:math.AG/0601395.

\bibitem{McKay}
John McKay.
\newblock Graphs, singularities, and finite groups.
\newblock In {\em The Santa Cruz Conference on Finite Groups (Univ. California,
  Santa Cruz, Calif., 1979)}, volume~37 of {\em Proc. Sympos. Pure Math.},
  pages 183--186. Amer. Math. Soc., Providence, R.I., 1980.

\bibitem{Morrison}
David~R. Morrison.
\newblock The birational geometry of surfaces with rational double points.
\newblock {\em Math. Ann.}, 271(3):415--438, 1985.

\bibitem{Nakamura}
Iku Nakamura.
\newblock Hilbert schemes of abelian group orbits.
\newblock {\em J. Algebraic Geom.}, 10(4):757--779, 2001.

\bibitem{Pandharipande-Thomas1}
Rahul Pandharipande and Richard Thomas.
\newblock {Curve counting via stable pairs in the derived category}.
\newblock arXiv:math/0707.2348.

\bibitem{Johnson-Panharipande-Tseng}
Hsian-Hua~Tseng Paul~Johnson, Rahul~Pandharipande.
\newblock {Abelian Hurwitz-Hodge integrals}.
\newblock arXiv:math/0803.0499.

\bibitem{Reid-asterisque}
Miles Reid.
\newblock La correspondance de {M}c{K}ay.
\newblock {\em Ast\'erisque}, (276):53--72, 2002.
\newblock S\'eminaire Bourbaki, Vol.\ 1999/2000.

\bibitem{Ruan-crepant}
Yongbin Ruan.
\newblock The cohomology ring of crepant resolutions of orbifolds.
\newblock In {\em Gromov-Witten theory of spin curves and orbifolds}, volume
  403 of {\em Contemp. Math.}, pages 117--126. Amer. Math. Soc., Providence,
  RI, 2006.

\bibitem{Zhou-An}
Jian Zhou.
\newblock {Crepant resolution conjecture in all genera for type A
  singularities}.
\newblock arXiv:0811.2023.

\bibitem{Zhou-HH}
Jian Zhou.
\newblock {On computations of Hurwitz-Hodge integrals}.
\newblock arXiv:math/0710.1679.

\end{thebibliography}
\bibliographystyle{plain}
\end{document}